\documentclass[12pt]{amsart}
\PassOptionsToPackage{usenames,dvipsnames}{xcolor}
\usepackage{txfonts,fancyvrb}
\usepackage{pgf,pgfplots}
\pgfplotsset{compat=1.14}
\usepackage{mathrsfs}
\usepackage{algorithm}
\usepackage[noend]{algpseudocode}
\makeatletter
\def\BState{\State\hskip-\ALG@thistlm}
\makeatother

\usepackage{subfig}
\usepackage{amscd,amssymb,mathabx}
\usepackage[arrow,matrix,graph,frame,poly,arc,tips]{xy}
\usepackage{graphicx,calrsfs}
\usepackage{mathtools}

\usepackage{mdwlist}
\usepackage{multirow}
\usepackage{enumerate}
\usepackage{verbatim}
\usepackage{etoolbox}
\AtEndEnvironment{proof}{\setcounter{claim}{0}}

\usepackage{pinlabel}

\DeclarePairedDelimiter{\form}{\langle}{\rangle}

\newcommand\be{\begin{enumerate}}
\newcommand\ee{\end{enumerate}}
\newcommand\bp{\begin{proof}}
\newcommand\ep{\end{proof}}
\newcommand\bpp{\begin{prop}}
\newcommand\epp{\end{prop}}
\newcommand\bpb{\begin{prob}}
\newcommand\epb{\end{prob}}
\newcommand\bd{\begin{defn}}
\newcommand\ed{\end{defn}}
\newcommand\bh{\begin{hint}}
\newcommand\eh{\end{hint}}

\newcommand\fform[1]{\langle\!\langle #1\rangle\!\rangle}

\newcommand\fC{\mathfrak{C}}
\newcommand\fT{\mathfrak{T}}
\newcommand\bC{\mathbb{C}}

\newcommand\bR{\mathbb{R}}

\newcommand\bQ{\mathbb{Q}}
\newcommand\bZ{\mathbb{Z}}
\newcommand\bH{\mathbb{H}}

\renewcommand\Im{\operatorname{Im}}

\newcommand\Isom{\operatorname{Isom}}

\newcommand\Id{\operatorname{Id}}
\newcommand\SL{\operatorname{SL}}

\newcommand\PSL{\operatorname{PSL}}

\newcommand\sse{\subseteq}
\newcommand\co{\colon}
\DeclareMathOperator\tr{tr}

\newcommand\Rot{\operatorname{Rot}}

\def\thetitle{Shapes of hyperbolic triangles and once-punctured torus groups}
\def\theauthors{{Sang-hyun Kim, Thomas Koberda, Jaejeong Lee, Ken'ichi Ohshika, and Ser Peow Tan}}
\usepackage{hyperref}
\hypersetup{
   colorlinks=false,
   plainpages,
   urlcolor=black,
   linkcolor=black
   pdftitle=   \thetitle,
   pdfauthor=  {\theauthors}
}

\theoremstyle{plain}
\newtheorem{thm}{Theorem}[section]
\newtheorem{lem}[thm]{Lemma}
\newtheorem{cor}[thm]{Corollary}
\newtheorem{prop}[thm]{Proposition}

\newtheorem{que}[thm]{Question}
\newtheorem*{claim*}{Claim}

\newtheorem{queAlph}{Question}

\theoremstyle{remark}

\newtheorem{rem}[thm]{Remark}

\theoremstyle{definition}
\newtheorem{defn}[thm]{Definition}
\newtheorem{prob}{Problem}[section]

\begin{document}
\title\thetitle



\author[S. Kim]{Sang-hyun Kim}
\address{School of Mathematics, Korea Institute for Advanced Study (KIAS), Seoul, 02455, Korea}
\email{skim.math@gmail.com}
\urladdr{https://www.cayley.kr}

\author[T. Koberda]{Thomas Koberda}
\address{Department of Mathematics, University of Virginia, Charlottesville, VA 22904-4137, USA}
\email{thomas.koberda@gmail.com}
\urladdr{http://faculty.virginia.edu/Koberda}

\author[J. Lee]{Jaejeong Lee}
\address{Institute of Basic Science, Sungkyunkwan University, Suwon, 16419, Korea}
\email{j.lee@skku.edu}

\author[K. Ohshika]{Ken'ichi Ohshika}
\address{Department of Mathematics, Gakushuin University}
\email{ohshika@math.gakushuin.ac.jp}

\author[S. P. Tan]{Ser Peow Tan\\ {\tiny \it with an appendix by} Xinghua Gao}
\address{Department of Mathematics, National University of Singapore, Singapore}
\email{mattansp@nus.edu.sg}
\urladdr{www.math.nus.edu.sg/~mattansp}

\address{School of Mathematics, Korea Institute for Advanced Study (KIAS), Seoul, 02455, Korea}
\email{xgao29@mail.illinois.edu}
\maketitle

\begin{abstract}
Let $\Delta$ be a hyperbolic triangle with a fixed area $\varphi$.
We prove that for all but countably many $\varphi$, 
generic choices of $\Delta$ have the property that the group generated by  the $\pi$--rotations
about the midpoints of the sides of the triangle admits no nontrivial relations. 
By contrast, we show for all $\varphi\in(0,\pi)\setminus\bQ\pi$, a dense set of triangles does afford nontrivial relations,
 which in the generic case map to hyperbolic translations.
To  establish  this fact, we study the deformation space $\fC_\theta$ of singular hyperbolic metrics on a torus with a single cone point of angle $\theta=2(\pi-\varphi)$, and answer an analogous question for the holonomy map $\rho_\xi$ of such a hyperbolic structure $\xi$. 
In an appendix by X.~Gao, concrete examples of $\theta$ and $\xi\in\fC_\theta$ are given where the image of each
$\rho_\xi$ is finitely presented, non-free and torsion-free; in fact, those images will be isomorphic to the fundamental
groups of closed hyperbolic 3--manifolds.
\end{abstract}
\date\today

\maketitle

\section{Introduction}
Take a 
geodesic triangle in the hyperbolic plane, and consider the rotations of angle $\pi$ about the midpoints of the three sides,
which we call the \emph{side involutions}. It is natural to wonder whether or not some nontrivial compositions of the
side involutions will move the triangle exactly back to itself. 

For a regular pentagon in the Euclidean plane, there are many ``unexpected'' coincidences
of this type~\cite{Schwartz2002AMM}.
On the other hand, a regular tetrahedron in $\mathbb{R}^3$ can never be moved back to its original position~\cite{AD2018-tetrahedron}.

For a hyperbolic triangle the answer depends on, among other  things, its area. 
For instance, if the area is a rational multiple of $\pi$, then so is the interior angle sum. In this case, 
a suitable power of a composition of the three side involutions will be trivial. 
Even if the area is not a rational multiple of $\pi$, torsion relations still appear for a dense choice of hyperbolic
triangles (cf.~Theorem~\ref{t:torsion}  below). So, we are led to address the question whether or not relations that are not
``consequences'' of torsion relations can still be found for a dense set of  triangles;
we refer the reader to Definition~\ref{d:torsion-type} for a precise formulation.

\begin{queAlph}\label{q:main1}
In the space of hyperbolic triangles whose area is fixed as $\varphi\in(0,\pi)$,
what are the necessary and sufficient conditions for the side involutions to generate the Coxeter group
$\mathbb{Z}_2\ast\mathbb{Z}_2\ast\mathbb{Z}_2$? 
Do the involutions for a dense set of hyperbolic triangles admit relations that are not consequences of torsion relations?
\end{queAlph}
Here, the space of hyperbolic triangles $ABC$ with fixed area $\varphi$ is identified with its parameter space of the interior angles:
\[\fT_\varphi:=\{(\theta_A,\theta_B,\theta_C)\in(0,\pi)^3 \mid \theta_A+\theta_B+\theta_C=\pi-\varphi\}.\]

Our take on Question~\ref{q:main1} will be via the space of marked incomplete hyperbolic structures on a punctured torus.
{Throughout this paper}, we let $M$ be a torus minus a puncture $p$.
Let us fix an angle $\theta\in(0,2\pi)$. We consider the deformation space $\fC_\theta$, called  a \emph{Fricke--Klein space}~\cite{Ji2016},
of marked incomplete hyperbolic structures on $M$ having exactly one conical singularity at the puncture $p$ with a cone angle $\theta$. 
We will also consider the Fricke--Klein space  \[\fC:=\bigcup_{\theta\in(0,2\pi)}\fC_\theta\]
of all marked hyperbolic structures with cone angles in $(0,2\pi)$ on $M$.

We fix an oriented meridian $X$ and an oriented longitude $Y$ of $M$, and identify $\pi_1(M)$ with a rank-two free group
\[ F_2:=\form{X,Y,Z\mid XYZ=1}.\] 
It is convenient to regard $F_2$ as an index-two subgroup of a free Coxeter group
\[
W:=\bZ_2\ast\bZ_2\ast\bZ_2=\form{P,Q,R\mid P^2=Q^2=R^2=1}\]
using the embedding
\[X\mapsto QR,\ Y\mapsto RP,\ Z\mapsto PQ.\]

Each hyperbolic structure $\xi\in\fC_\theta$ defines a conjugacy class of \emph{holonomy maps}, a representative of which we denote by
\[\rho_\xi\co \pi_1(M)\to\PSL(2,\bR)\cong\Isom_+(\bH^2).\] These holonomy maps are defined via the developing map $\tilde M\to\bH^2$.
The isometries $\{\rho_\xi(X),\rho_\xi(Y),\rho_\xi(Z)\}$ are hyperbolic elements with pairwise intersecting axes. 
Because the group $\Im\rho_\xi$ contains two hyperbolic elements with intersecting axes, it is always non-elementary. 
Since $F_2$  is Hopfian, we have that $\rho_\xi$ is faithful if and only if   $\Im\rho_\xi$ nonabelian and free. 

As is well-known (see e.g.\cite{Goldman2003GT}), the space $\fC_\theta$ can be identified with its parameter space
\[\{(x,y,z)\in(2,\infty)^3\mid x^2+y^2+z^2-xyz-2= -2\cos (\theta/2)\}\]
by the \emph{character} map
\[\xi\mapsto \left(\tr^+ \rho_\xi(X),\tr^+ \rho_\xi(Y),\tr^+ \rho_\xi(Z)\right).\]
Here, $\tr^+ (\cdot)$ is understood as the absolute value of the trace of an element in $\PSL(2,\bR)$.

\begin{figure}[htb!]
	\labellist
	\pinlabel {\tiny $\gamma_X$} at 380 210
	\pinlabel {\tiny $\gamma_Y$} at 150 374
	\pinlabel {\tiny $\gamma_Z$} at 116 60
	\pinlabel {\tiny $h_X$} at 239 230
	\pinlabel {\tiny $O_P$} at 220 124
	\pinlabel {\tiny $O_Q$} at 306 197
	\pinlabel {\tiny $O_R$} at 155 225
	\pinlabel {\tiny $\delta_X$} at 260 154
	\pinlabel {\tiny $O_A$} at 295 355
	\pinlabel {\tiny $O_B$} at 80 85
	\pinlabel {\tiny $O_C$} at 340 95
	\pinlabel {\tiny $R(O_C)$} at 66 290
	\pinlabel {\tiny $\theta_A$} at 266 280
	\pinlabel {\tiny $\theta_B$} at 246 296
	\pinlabel {\tiny $\theta_C$} at 296 270
	\endlabellist
	\centering
	\includegraphics[width=0.4\textwidth]{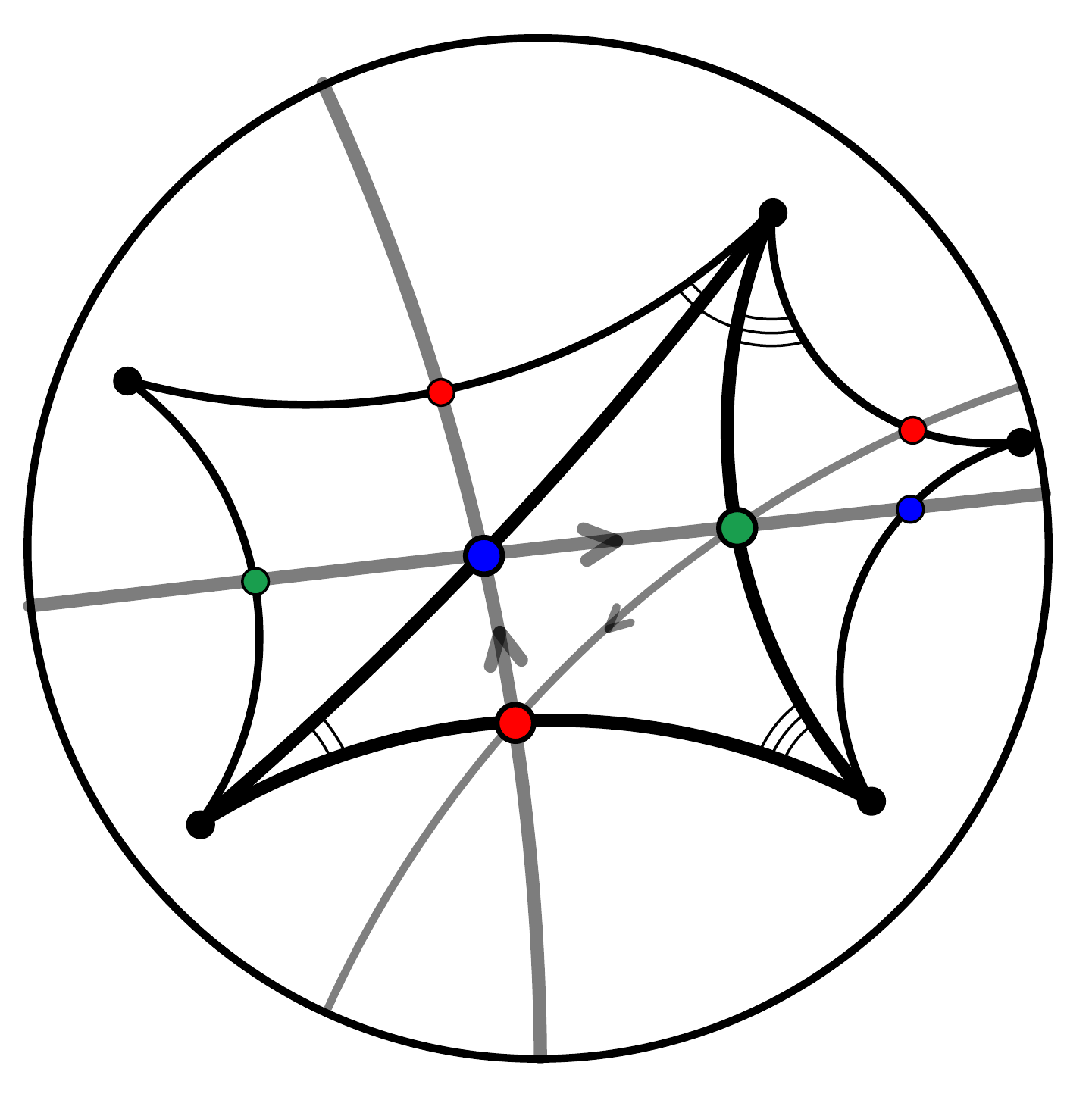}
	\caption{Developing $M$ to $\bH^2$.}
	\label{f:setup-small}
\end{figure}

To see the connection between $\fC_\theta$ and Question~\ref{q:main1},
pick a hyperbolic structure $\xi\in\fC_\theta$ on $M$ for some fixed $\theta\in(0,2\pi)$. 
Let $\gamma_X$ be the simple closed geodesic realizing $X\in \pi_1(M)$. There exists a unique incomplete geodesic $\delta_X$ on
$M$ disjoint from $\gamma_X$, and that starts and ends at $p$. Construct the geodesics 
$\{\gamma_Y,\gamma_Z,\delta_Y,\delta_Z\}$ analogously.
By cutting $M$ along $\delta_X\cup\delta_Y\cup\delta_Z$, we obtain two isometric hyperbolic triangles whose
interior angle sums are $\theta/2$. The shape of this triangle determines a point in $\fT_{\pi-\theta/2}$.
A lift of this picture to $\bH^2$ is illustrated in  Figure~\ref{f:setup-small}.

As in the figure, we let $\{O_P, O_Q, O_R\}$ denote the midpoints of the sides from one of the triangles.
We write $\{\hat\rho_\xi(P),\hat\rho_\xi(Q),\hat\rho_\xi(R)\}$ for the corresponding side involutions.
Then we have that  
\[\rho_\xi(X)=\hat\rho_\xi(Q)\cdot \hat\rho_\xi(R),\quad  \rho_\xi(Y)=\hat\rho_\xi(R)\cdot \hat\rho_\xi(P),\quad
\rho_\xi(Z)=\hat\rho_\xi(P)\cdot \hat\rho_\xi(Q).\]
Thus $\rho_\xi$ extends to a surjection
\[
\hat\rho_\xi\co W\to \form{\hat\rho_\xi(P),\hat\rho_\xi(Q),\hat\rho_\xi(R)},\]
called a \emph{Coxeter extention} of $\rho_\xi$, which is determined  uniquely from the holonomy representation
$\rho_\xi$~\cite{Goldman2003GT}. 
Moreover, the spaces $\fT_{\pi-\theta/2}$ and  $\fC_\theta$ are homeomorphic
via a homeomorphism which conjugates this extension (Lemma~\ref{l:parameters}).

Through Coxeter extensions, we will approach Question~\ref{q:main1} by studying the following closely related
question.
\begin{queAlph}\label{q:main2}
For a fixed angle $\theta\in(0,2\pi)$,
what are the necessary and sufficient conditions for $\xi\in\fC_\theta$ to yield a faithful representation $\rho_\xi\co \pi_1(M)\to\PSL(2,\bR)$? 
Is there a dense choice of $\xi\in\fC_\theta$ for which the representation $\rho_\xi$ admits 
relations that are not consequences of torsion relations?
\end{queAlph}

Let $A$ be a topological space. We say that a certain statement holds for a \emph{very general} (or, \emph{generic}) $\xi\in A$ if
 the statement holds for all $\xi$ in some dense $G_\delta$ subset of $A$. If $A$ is an open subset of a Euclidean space,
 then it follows that the statement holds on a full-measure subset of $A$.

In this paper, we are primarily interested in establishing the faithfulness part of Question~\ref{q:main2}
for a very general point $\xi\in\fC_\theta$, and in answering the density question for a suitable interpretation of the phrase
``consequences of torsion relations''.

\begin{rem}
One may also consider \emph{side reflections} instead of side involutions of a hyperbolic triangle $\Delta$ with angles  $(\theta_A,\theta_B,\theta_C)$, although we will not pursue this direction. We briefly note in this case that two representations $\hat\rho_\xi$ and $\rho_\xi$ are still well-defined on $ W$ and on $F_2$ respectively.
The representation $\rho_\xi$ is a holonomy map of an incomplete hyperbolic three-punctured sphere $S_{0,3}$ with three conical singularities of angles $(2\theta_A,2\theta_B,2\theta_C)$.
This approach could be useful when one is interested in periodic billiard orbits on $\Delta$. 
\end{rem}

The faithfulness question is easier when the cone angle $\theta$ is not fixed. Namely, we have that  $\rho_\xi$ is faithful for a very general point $\xi\in \fC$ (cf. Proposition~\ref{p:free}). So, we will be actually interested in the case when $\theta$ is fixed.

Our first main result answers the first half of Question~\ref{q:main2} with probability one:

\begin{thm}[cf. Theorem~\ref{t:free}]\label{t:free-intro}
If $\theta\in(0,2\pi)$ has  the property that $\cos\theta$ is transcendental, then  
a very general point $\xi$ in $\fC_\theta$ corresponds to a faithful holonomy map. 
\end{thm}
Note that the hypothesis holds for all but countably many values of $\theta$.  The reader may compare Theorem~\ref{t:free-intro} to other generic phenomena in $\PSL(2,\bR)$; see~\cite{DK2006Duke,FW2007,KKM2019}, the first two of which deal with the faithfulness question in the case of closed surface groups.
One can trace at least back to Hausdorff~\cite{Hausdorff1949} the idea of using the transcendency of $\cos\theta$ in order to produce faithful group actions.

It is a well-known consequence of the Gelfond--Schneider theorem that $\cos\theta$ is transcendental if $\theta/\pi$ is irrational and algebraic. Thus, we deduce more concrete cases when $\rho_\xi$ is faithful as follows.
\begin{cor}\label{c:algirr}
If $\theta\in(0,2\pi)$ is an irrational algebraic multiple of $\pi$, 
then $\rho_\xi$ is faithful for a very general point $\xi\in \fC_\theta$.
\end{cor}

We have briefly mentioned the abundance of torsion relations at the beginning of this introduction. Indeed,
the image of $\rho_\xi$ is dense in $\PSL(2,\bR)$ unless $\rho_\xi$ is discrete (and hence Fuchsian), as is true of all Zariski dense subgroups
of simple Lie groups~\cite{Katok1992book}. 
Since elliptic isometries form an open subset in $\PSL(2,\bR)$, it is reasonable to expect that an elliptic element in the image of
$\rho_\xi$ would not have constant rotation angle under perturbations in $\fC_\theta$ near $\xi$. Thus, a torsion relation would
appear under small deformations of $\xi$ having a fixed cone angle $\theta$; we direct the reader to Theorem~\ref{t:torsion}
for a concrete formulation of this intuition. 
Moreover, in the presence of torsion there are generally many more elements of $\ker\rho_\xi$ that are consequences of such torsion relations. 

We aim to find relations of $\rho_\xi$ which are not ``consequences'' of torsion relations. 
To state the result rigorously, we make the following definition. 
\bd\label{d:torsion-type}
A word $u\in F_2$ is of \emph{torsion-type}
if $u$ belongs to the set
\[
\bigcup_{m\ge2}\bigcup_{r\in F_2\setminus\{1\}} \fform{r^m}.\]
A word is of \emph{non-torsion-type} if it is not of torsion-type.
\ed

\

For a word $u\in F_2$, being a non-torsion-type word is a stronger condition than not equal to a proper power.  
In the special case when $\theta=2\pi/n$, the group $\rho_\xi(F_2)$ is isomorphic to a Fuchsian orbifold group
\[\form{X,Y\mid [X,Y]^n=1}\]
for all choices of $\xi\in\fC_\theta$; in this case, every kernel element must be of torsion-type. 
When   $\theta\not\in\bQ\pi$, we can produce non-torsion-type kernel elements as addressed in the second half of Question~\ref{q:main2}. 

\begin{thm}[cf.~Theorem~\ref{t:nontorsion}]\label{t:nontorsion-intro}
If $\theta\in(0,2\pi)$ is an irrational multiple of $\pi$,
then $\ker\rho_\xi<F_2$ contains a palindromic non-torsion-type word for a dense choice of $\xi\in\fC_\theta$.
\end{thm}

In particular, when $\ker\rho_\xi$ contains a non-torsion-type word, 
the image $\rho_\xi(F_2)$ cannot be isomorphic to a one-relator group with torsion with respect to the
generating set $\{\rho_\xi(X),\rho_\xi(Y)\}$.

We will actually produce a dense subset $S_\theta\sse (2,\infty)$ such that each point on the coordinate curves 
\[ \{x=s\}, \quad \{y=s\},\quad \{z=s\}\] in $\fC_\theta$ for $s\in S_\theta$ corresponds to holonomy maps admitting 
non-torsion-type kernel elements. 
The points on the intersection of any two of these coordinate curves ({\it double points}) form a dense subset of $\fC_\theta$ with the same property.
\begin{rem}
A consequence of Theorem~\ref{t:nontorsion-intro} is that a free indiscrete representation is
\emph{algebraically unstable}, even when restricted to $\fC_\theta$; roughly speaking this means that such a representation
is a limit of non-free representations.  In the absence of such a ``relativizing'' restriction, this type of instability
is well-known for non-abelian free groups, even in a much more general setting of
connected Lie groups~\cite{Glutsyuk2011TG}. On the other hand, Sullivan~\cite{Sullivan1985QC2}
established algebraic stability of all convex--cocompact subgroups in  $\PSL(2,\bC)$.
\end{rem}

Let us direct our attention back to Question~\ref{q:main1}. 
The faithfulness part has the same answer as Question~\ref{q:main2}; see~Proposition~\ref{p:equiv}(1).
Regarding the density part, it is unclear to the authors whether the non-torsion-type kernel elements in $F_2$ found in
 Theorem~\ref{t:nontorsion-intro} can still be non-torsion-type in $W$. However, their additional property of
 being palindromic has a simple interpretation in $W$: they are products of two involutions in $W$ (Lemma~\ref{l:pal}).
Thus we establish the following partial answer to Question~\ref{q:main1}.

\begin{cor}\label{cor:main1}
The following conclusions hold:
\be
\item
For all but countably many $\varphi\in(0,\pi)$, the side involutions of a very general triangle with area $\varphi$
generates the Coxeter group  $W$. 
\item
If $\varphi$ is an irrational multiple of $\pi$, then for a dense choice of hyperbolic triangles with area $\varphi$
the side involutions admit a relation which is the product of two involutions in $W$. \ee
\end{cor}

We remark that if $u\in W$ is a product of two involutions then, for all $\xi\in\fC$, the image $\hat\rho_\xi(u)$
is either trivial or hyperbolic (Proposition~\ref{p:equiv}(2)). In fact, it is hyperbolic for a very general point in
$\fC_\theta$ (Remark~\ref{rem:palindrome}).

As far as the authors are aware, no examples of $\xi\in\fC$ were
previously documented such that $\rho_\xi(F_2)$ is non-free and torsion-free.

\begin{que}\label{q:main3}
If $\theta\in(0,2\pi)$ is an irrational multiple of $\pi$,  under what conditions on $\xi$ is $\rho_\xi$ non-faithful with torsion-free image?
\end{que}
In the appendix written by X.~Gao, computational heuristics for approaching Question~\ref{q:main3}, together with their implementation,
 are given.
The examples of pairs $(\theta,\xi)$ exhibited therein have the property that  the images of $\rho_\xi$ are actually
isomorphic to fundamental groups of closed hyperbolic 3--manifolds.

\section{Fricke--Klein space}\label{s:prelim}
Most of the material in the section is well-known; we direct the reader to~\cite{Goldman2003GT} for a standard reference.
We adopt the standing convention that group and matrix actions are on the left, unless stated otherwise.

The group $\SL(2,\bR)$ acts on $\bH^2$ by M\"obius transforms, with kernel given by the center, and with image
\[\PSL(2,\bR)\cong\Isom_+(\bH^2).\]
Let $\tr(A)$ denote the trace of a matrix $A$. For each word $w=w(X,Y)\in F_2=\form{X,Y}$ there exists a \emph{trace polynomial}
\[g_w\in\bZ[x,y,z]\] such that whenever $U, V\in\SL(2,\bR)$ we have
\[\tr w(U, V) =g_w(\tr U,\tr V,\tr UV).\] The existence of this polynomial is one of the simplest instances of invariant
theory on character varieties~\cite{Procesi76}; see also~\cite{Jorgensen1982} for a concrete formula which computes $g_w$. 

For example, the trace polynomial of $[X,Y]$ is easily seen to be
\[
g_{[X,Y]}(x,y,z)=\kappa(x,y,z):=x^2+y^2+z^2-xyz-2\in\bZ[x,y,z].\]
 \begin{lem}\label{l:nonconst}
For each $w\in F_2\setminus\{1\}$, the polynomial $g_w(x,y,z)$ is not constant.
\end{lem}
\bp
Let $w=w(X,Y)\in F_2=\form{X,Y}$.
If $g_w$ were constant, we would have \[g_w(x,y,z)=g_w(2,2,2)=\tr w(\Id,\Id)=\tr \Id =2.\] 
However, there exist two--generated nonelementary purely loxodromic Fuchsian groups which are nonabelian and free, namely, Schottky groups of rank two. In particular,
there are choices of parameters defining $X$ and $Y$ such that the element $w$ is loxodromic and therefore has trace different from
$2$.
This immediately implies  $g_w$ cannot be identically  equal to $2$.\ep

In this paper, we are mostly concerned with matrices in $\PSL(2,\bR)$. The trace of such a matrix $A\in\PSL(2,\bR)$
is only determined up to sign, so as mentioned in the introduction, we often use the quantity
\[\tr^+ A := |\tr A|.\]

Recall we let $\fC_\theta$ denote the deformation space (\emph{Fricke--Klein space}) of marked incomplete hyperbolic metrics on
$M$ with a fixed cone angle $\theta\in(0,2\pi)$ at $p$.
Each point $\xi\in\fC_\theta$ corresponds to a conjugacy class of a holonomy map
\[\rho\co \pi_1(M)=\form{X,Y,Z\mid XYZ=1}\to { \PSL(2,\bR)}\]
so that $\rho[X,Y]$ is a rotation of angle $\theta$.
The chosen representation $\rho$ has a unique lift to
\begin{align}\label{e:lift}
\tilde \rho\co \pi_1(M)\to \SL(2,\bR)
\end{align}
which satisfies \[\tr\tilde\rho(X),\tr\tilde\rho(Y),\tr\tilde\rho(Z)>2.\]
It turns out~\cite{Goldman2003GT} that the lift $\tilde\rho$ satisfies
\[\tr\tilde\rho[X,Y] = -2\cos(\theta/2).\]
Using the parametrization
\[
\xi\mapsto (\tr\tilde\rho(X),\tr\tilde\rho(Y), \tr\tilde\rho(Z))
=  (\tr^+ \rho(X),\tr^+ \rho(Y), \tr^+ \rho(Z)),\]
and the results of Section 3 of \cite{Goldman2003GT}, we can therefore identify
\[\fC_\theta=\{(x,y,z)\in(2,\infty)^3\mid \kappa(x,y,z)= -2\cos (\theta/2)\}.\]

\begin{rem}\label{r:normal}
The surjectivity of the above parametrization can be seen by defining a \emph{normal form}
$\tilde\rho_\xi$ of $\xi=(x,y,z)\in(2,\infty)^3$. We let $\zeta\ge1$ be the unique real number satisfying $\zeta+\zeta^{-1}=z$
and define $\tilde\rho_\xi$ via
\[\tilde\rho_\xi(X):=\begin{pmatrix} x & -1 \\ 1 & 0\end{pmatrix},\quad
\tilde\rho_\xi(Y):=\begin{pmatrix} 0 & \zeta^{-1} \\ { -\zeta} & y\end{pmatrix}.\]
From an easy computation~\cite[Section 2.2.3]{Goldman2009} one readily sees that
\[\tr\tilde\rho_\xi(Z)=\tr\rho_\xi((XY)^{-1})=z.\]
We let $\rho_\xi\co \pi_1(M)\to\PSL(2,\bR)$ be the projection of $\tilde\rho_\xi$.
These concrete formula for $\tilde\rho_\xi$ and $\rho_\xi$ are not actually needed in this paper.
\end{rem}

\section{Very general representations are free}

Recall that
\[
\fC:=\bigcup_{\theta\in(0,2\pi)}\fC_\theta=\{(x,y,z)\in(2,\infty)^3\mid \kappa(x,y,z)\in(-2,2)\}.\]
 The following is well-known to experts;  we include a proof as motivation for the next theorem.
\begin{prop}\label{p:free}
For a very general point $\xi\in \fC$, the representation $\rho_\xi$ is faithful.
\end{prop}

\bp
For each $w\in F_2$, let us set
\[Y_w:=\{\xi\in \fC\mid \rho_\xi(w)=\Id\}\sse\{\xi\in\fC\mid g_w(\xi)=\pm2\}.\]
By~Lemma~\ref{l:nonconst}, we see that $Y_w$ is contained in a proper algebraic subset of $\bR^3$ whenever
$w\ne1$. In particular, $Y_w$ has no interior inside the open subset $\fC$ of $\bR^3$. It follows that
\[\fC\setminus\left(\bigcup_{w\in F_2\setminus\{1\}}Y_w\right)\] is a $G_\delta$--set.
\ep
\begin{rem}
The above proof actually implies that a very general point $\xi\in \fC$ corresponds to representations
$\rho_\xi$ without any nontrivial parabolic elements in the image. 
\end{rem}

\begin{thm}\label{t:free}
If $\theta\in(0,2\pi)$ has the property that $\cos\theta$ is transcendental, 
then for a very general point $\xi$ in $\fC_\theta$ 
the representation $\rho_\xi$ is faithful.
\end{thm}

\begin{rem}
Note that the conclusion of Theorem~\ref{t:free} holds for all but countably many $\theta$.
Observe furthermore  that the  theorem does not extend to all $\theta\in(0,2\pi)$.
For instance, if $\theta\in(\bQ\setminus\bZ)\pi$ then $\rho_\xi(F_2)$ has nontrivial torsion.
\end{rem}

We let $\bar\bQ$ denote the algebraic closure of $\bQ$ in $\bC$.

\bp[Proof of Theorem~\ref{t:free}]
Let $\theta$ be given as in the hypothesis. Recall for $w\in F_2$ we defined the trace polynomial $g_w(x,y,z)$.
The following claim is the key  observation in the proof.
\begin{claim*}
For each $w\in F_2\setminus\{1\}$ and for each algebraic number $c$, 
the set 
\[ Y_w(c):=\{\xi\in \fC_\theta\mid g_w(\xi)=c\}\]
has empty interior in $\fC_\theta$.
\end{claim*}

Let us assume the claim for now. 
Then each $\xi$ in the dense $G_\delta$-set
\[
\fC_\theta\setminus\left(\bigcup_{w\in F_2\setminus\{1\}} Y_w(\pm2)\right)\]
has the property that $\rho_\xi$ is injective, which completes the proof.

We now establish the claim. Write
\[
g_w(x,y,z) = \sum_{i=0}^n f_i(x,y)z^i\]
for suitable $f_i\in\bZ[x,y]$. We first consider the case $n\ge1$ and $f_n\ne0$.
We define
\[ S:=\{(x,y,z)\in \fC_\theta\cap(\bQ\times\bQ\times\bR)\mid f_n(x,y)\ne0\}.\]
Since $\bQ^2$ is dense in $\bR^2$, the set $S$ is dense in $\fC_\theta$.
For all $\xi=(x,y,z)\in S$, we see that $z\not\in \bar\bQ$ from the transcendency hypothesis on $\cos\theta$. 
So, we have $g_w(\xi)\not\in\bar\bQ$.
In particular, we obtain $S\cap Y_w(c)=\varnothing$. Since $S$ is dense, the claim is proved in this case.

We now consider the case that $n=0$. We have $g_w=f_0(x,y)$. By Lemma~\ref{l:nonconst}, we know that $f_0$ is not constant.
Mimicking the previous argument, we consider the following two dense subsets of $\fC_\theta$: 
\begin{align*}
S_1&:=\{(x,y,z)\in \fC_\theta\cap(\bQ\times\bR\times\bQ)\mid f_0(x,y)\ne0\},\\
S_2&:=\{(x,y,z)\in \fC_\theta\cap(\bR\times\bQ\times\bQ)\mid f_0(x,y)\ne0\}.
\end{align*}
Again by considering transcendency of coordinates, we obtain that
$g_w(S_i)\cap\bar\bQ=\varnothing$ for $i=1$ or $i=2$; that is, we have  
$S_i\cap Y=\varnothing$. The conclusion follows from the density of $S_i$.
\ep
\begin{rem}
In Proof of  Theorem~\ref{t:free}, we need Claim only for  the case $c=\pm 2$.
Since $g_w$ is an analytic function, to show that $Y_w(c)$ has empty interior for general $c$, it is sufficient to show that $g_w$ is not constant in $\fC_\theta$.
We chose an argument different from this, imposing an extra condition that $c$ is algebraic,  to give specific points  not contained in $Y_w(c)$, with a concrete description.
\end{rem}

If $\theta$ is a rational multiple of $\pi$, then the image of $\rho_\xi$ contains nontrivial torsion elements, so that
$\Im\rho_\xi$ is non-free regardless of the choice of $\xi\in \fC_\theta$.
On the other hand,  if $\theta$ is an algebraic irrational multiple of $\pi$ then the Gelfond--Schneider theorem implies that
$\cos\theta$ is transcendental, so that $\rho_\xi$ is faithful for a very general point $\xi\in\fC_\theta$ (Theorem~\ref{t:free}).
Hence, the only remaining case in the faithfulness part of Question~\ref{q:main2} is the following.

\begin{que}\label{q:trans}
If $\theta\in(0,2\pi)$ is a transcendental multiple of $\pi$, 
and if $\cos\theta$ is algebraic,
then is $\rho_\xi$  faithful for a very general point $\xi\in \fC_\theta$?
\end{que}

The general case of this question seems mysterious. For instance, one may ask whether or not a very general point $\xi$ in
$\fC_{\arccos(1/3)}$ correspond to a faithful holonomy map. 

\section{Coxeter extensions}\label{s:coxeter}
It will be computationally convenient for us to consider an embedding from $\rho_\xi(F_2)$ to a
bigger group generated by involutions (as in Question~\ref{q:main2}). 
Such an embedding will be also used here to see the connection between Questions A and B.

Let us begin with an algebraic discussion. We say a word $w(X,Y)$ is palindromic in $\{X,Y\}$
if it reads the same forward and backward, that is, $w(X,Y)^{-1}=w(X^{-1},Y^{-1})$. We say \[w\in F_2=\form{X,Y,Z\mid XYZ=1}\] is
\emph{palindromic} if it can be expressed as being palindromic in either $\{X,Y\}$, $\{Y,Z\}$ or $\{Z,X\}$.

Recall we are regarding $F_2=\form{X,Y,Z\mid XYZ=1}$ as an index-two subgroup of $W=\form{P,Q,R}$ using the embedding
\[X\mapsto QR,\ Y\mapsto RP,\ Z\mapsto PQ.\]
Note, in particular, that $RXR=X^{-1}$ and $RYR=Y^{-1}$. 

\begin{lem}\label{l:pal}
A word $u\in F_2$ is palindromic if and only if it is a product of two involutions in $W$, with one of the two being $P$, $Q$, or $R$.
\end{lem}
We remark that a similar statement can be found in~\cite{Jorgensen1978} for the case when the ambient group is $\SL(2,\bC)$.
\bp[Proof of Lemma~\ref{l:pal}]
The word $u$ is palindromic in $\{X,Y\}$
if and only if
\[u(X,Y)^{-1}=u(X^{-1},Y^{-1})=u(RXR,RYR)=Ru(X,Y)R,\] which is true
if and only if
$(Ru)^2=1$. This last expression holds
if and only if
$Ru=gIg^{-1}$ for some $I\in\{P,Q,R\}$ and $g\in W$, as follows from the characterization of torsion in right-angled Coxeter groups.
It follows that $u=R\cdot gIg^{-1}$.

Similarly, one shows that $u$ is palindromic in $\{Y,Z\}$ (resp. $\{Z,X\}$) if and only if $u=P\cdot gIg^{-1}$ (resp. $u=Q\cdot gIg^{-1}$) for some $I\in\{P,Q,R\}$ and $g\in W$.
\ep

Turning to the geometric side of the embedding $F_2<W$, let us fix $\theta\in(0,2\pi)$ and $\xi=(x,y,z)\in \fC_\theta$.
Since $\kappa(x,y,z)\ne2$, the holonomy map $\rho_\xi$ is irreducible~\cite[Proposition~2.3.1]{Goldman2009}. 
It follows from \cite[Theorem B and Theorem 3.2.2]{Goldman2009} that 
 $\rho_\xi$ uniquely extends to  the \emph{Coxeter extension} $\hat\rho_\xi$  of $\rho_\xi$ as shown in the commutative diagram below:
\[
\xymatrix{
	W=\form{P,Q,R} \ar[rd]^{\hat\rho_\xi}\ar@{-}[d]\\
	F_2=\form{X,Y}\ar[r]^{\rho_\xi}&  \PSL(2,\bR)}
\]
Explicitly, one can define the desired extension by projecting the following formula~\cite{Goldman2009} in $\SL(2,\bR)$ to $\PSL(2,\bR)$:
\[ \hat\rho_\xi(R)=\frac{\tilde\rho_\xi(X)\tilde\rho_\xi(Y)-\tilde\rho_\xi(Y)\tilde\rho_\xi(X)}{\sqrt{2-\kappa(x,y,z)}},\]
where $\tilde\rho_\xi$ is the lift in \eqref{e:lift}.
The formulae for $P$ and $Q$ are analogous.

From now on, we will often suppress the symbols $\rho_\xi$ and $\hat\rho_\xi$, and simply write $w\in F_2$ or
$w\in W$ for $\hat\rho_\xi(w)$ when the meaning is clear.

We let \[2h_X,2h_Y,2h_Z\] be the translation lengths of the hyperbolic isometries $X,Y,Z\in\PSL(2,\bR)$, so that
\begin{align}\label{e:trace}
x=\tr X = 2\cosh h_X,
\end{align}
and similarly for $y$ and  $z$.
In the introduction we noted that the punctured torus $M$ can be cut along geodesics and developed
to $\bH^2$ as illustrated in Figure~\ref{f:setup-small}.

Let us continue to use the notation from the same figure. In particular, we let $\{\theta_A,\theta_B, \theta_C\}$
be the interior angles of the triangle $O_AO_BO_C$. 
We let $d_A$ (resp. $d_B$ and $d_C$) be the length of the geodesic segment $O_BO_C$ (resp. $O_CO_A$ and $O_AO_B$).
The isometry $X$ is a hyperbolic translation along the oriented geodesic $O_RO_Q$, and similar statements hold for the pairs
$(Y,O_PO_R)$ and $(Z,O_QO_P)$.
Moreover, the three isometries
\[ A := QPR,\quad B := RQP,\quad C := PRQ\]
are (counterclockwise) rotations of angle $\theta/2$ with centers $O_A$, $O_B$ and $O_C$, respectively.

The midpoints of the sides of this triangle are centers of the $\pi$-rotations  $P,Q,R$.
The quadrilateral $R(O_C)O_BO_CO_A$ can be regarded as a ``pseudo-fundamental domain''
for the incomplete hyperbolic structure $\xi\in\fC_\theta$ on $M$.
We refer the reader to~\cite{Goldman2003GT,Goldman2009} for more details.
We have that $\theta/2=\theta_A+\theta_B+\theta_C$, and
\[
K=[X,Y]=(QPR)^2=A^2\]
is a rotation by an angle $\theta$. We also record the identities \[B^2=[Y,Z], \quad C^2=[Z,X].\] 
If we denote by $d(\cdot,\cdot)$ the distances in $\bH^2$, then  the three lengths \[\{d(O_A,O_RO_Q), d(O_B,O_RO_Q),d(O_C,O_RO_Q)\}\]
 are all equal.

\begin{lem}\label{l:parameters}
Let $\theta\in(0,2\pi)$. Then there exists a homeomorphism $\Phi$ from
$\fT_{\pi-\theta/2}$ to $\fC_\theta$ such that 
for each $\eta\in\fT_{\pi-\theta/2}$, the Coxeter extension of $\Phi(\eta)$ is precisely the representation \[W\to\PSL(2,\bR)\]
coming from the side involutions. 
 More precisely, we have
\[
\cos\theta_A = \frac{
yz (x^2+2)- (y^2+z^2)x-x^3}
{2\sqrt{(xyz-x^2-y^2)(xyz-x^2-z^2)}}.\]
\end{lem}
\bp
Given $(\theta_A,\theta_B,\theta_C)\in\fT_{\pi-\theta/2}$, it is straightforward to define $\Phi$ using the hyperbolic law of cosines:
first find $\{d_A,d_B,d_C\}$, then $\{h_X,h_Y,h_Z\}$ and then $\{x,y,z\}$ using \eqref{e:trace}.

Now, given $(x,y,z)\in\fC_\theta$, an explicit formula of the inverse $\Phi^{-1}$ can be found as follows; see Figure~\ref{f:setup-small}. If we cut $M$ along the images of the incomplete geodesic $\delta_X=O_BO_C$ and the complete geodesic $\gamma_X$ extending $O_QO_R$, then we obtain two cylinders with a single common corner point $O$ (the image of $O_B$ or $O_C$) and boundary lengths $2h_X$ and $d_A$. By cutting the cylinder further along the shortest path from $O$ to the opposite boundary, we obtain a quadrilateral with a $2$-fold symmetry, which can be divided into two copies of a tri-rectangle with a non-right angle $\theta/4$ and opposite edge lengths $h_X$ and $d_A/2$, so that the edge of length $d_A/2$ is adjacent to the non-right angle.
We remark briefly that another way to obtain these tri-rectangles would be using a twist to be defined below.
After a suitable $X$-twist we may assume that the angles $\measuredangle O_QO_RO_P$ and $\measuredangle O_RO_QO_C$ are both right-angles, so that $O_CO_QO_RO_P$ is a tri-rectangle.

In this case we have a formula 
\[x/2=\cosh h_X=\cosh(d_A/2)\sin(\theta/4).\] See~\cite{Buser1992}, for example.
Once $\{d_A,d_B,d_C\}$ is found, using the hyperbolic law of cosines we can compute
\[\cos\theta_A = \frac{\cosh d_B\cosh d_C - \cosh d_A}{\sinh d_B \sinh d_C}.\] The remainder of the proof is straightforward.
\ep

We now exhibit a connection between the main questions in the introduction.
\begin{prop}\label{p:equiv}
For $\theta\in(0,2\pi)$ and $\xi\in\fC_\theta$, we have the following.
\be
\item
The representation $\rho_\xi$ is faithful if and only if  $\hat\rho_\xi$ is faithful.
\item
If $u\in F_2$ is palindromic, then the image $\hat\rho_\xi(u)\in\PSL(2,\bR)$ is either trivial or hyperbolic.
\ee
\end{prop}
\bp
(1) Suppose $\rho_\xi$ is faithful. If $1\ne g\in \ker\hat\rho_\xi$, then we have that
\[
g^2\in F_2\cap \ker\hat\rho_\xi=\ker\rho_\xi=\{1\}.\] Hence $g\in W $ has order two; this implies that $g$ is conjugate to
$P, Q$ or $R$. Without loss of generality, we may assume $g=P$. Then we have a contradiction since  $\hat\rho_\xi(P)=1$ and thus
\[\rho_\xi(Y^2)=\rho_\xi (RPRP)=\hat\rho_\xi (R^2)=1.\]
The converse is obvious.

(2) By Lemma~\ref{l:pal} if $u\in F_2$ is palindromic, then $u=I_1I_2$ for some involutions $I_1$ and $I_2$ in $W$.
Since the only nontrivial involutions in $\PSL(2,\bR)$ are $\pi$-rotations, so are the images $\hat\rho_\xi(I_1)$ and $\hat\rho_\xi(I_2)$.
Thus $\hat\rho_\xi(u)=\hat\rho_\xi(I_1)\hat\rho_\xi(I_2)$ is trivial (if $\hat\rho_\xi(I_1)=\hat\rho_\xi(I_2)$) or hyperbolic
(if $\hat\rho_\xi(I_1)\ne\hat\rho_\xi(I_2)$). 
\ep

\section{Density of non-faithful representations}

In this section we give geometric constructions of non-faithful holonomy maps.
The relations will be torsion elements as well as non-torsion-type.
The idea of construction is to consider the Coxeter extension of $F_2$ and deform a given representation $\rho$ slightly by \emph{twist}.

\subsection{Twist deformations}
Let us describe a deformation of a given structure $\xi\in\fC_\theta$ along one of the three coordinate curves
 \[
 \{x=\textup{const.}\}, \quad \{y=\textup{const.}\} , \quad \{z=\textup{const.}\}\]
passing through $\xi$ in $\fC_\theta$.

\begin{figure}[ht]
\labellist
\pinlabel {\tiny $O_P$} at 204 134
\pinlabel {\tiny $O_Q$} at 310 255
\pinlabel {\tiny $O_R$} at 195 255
\pinlabel {\tiny $O_A$} at 300 345
\pinlabel {\tiny $O_B$} at 70 120
\pinlabel {\tiny $O_C$} at 343 120
\pinlabel {\tiny $O_Q'$} at 236 195
\pinlabel {\tiny $O_R'$} at 135 255
\pinlabel {\tiny $O_A'$} at 180 365
\pinlabel {\tiny $L$} at 100 340
\endlabellist
\centering
\includegraphics[width=0.3\textwidth]{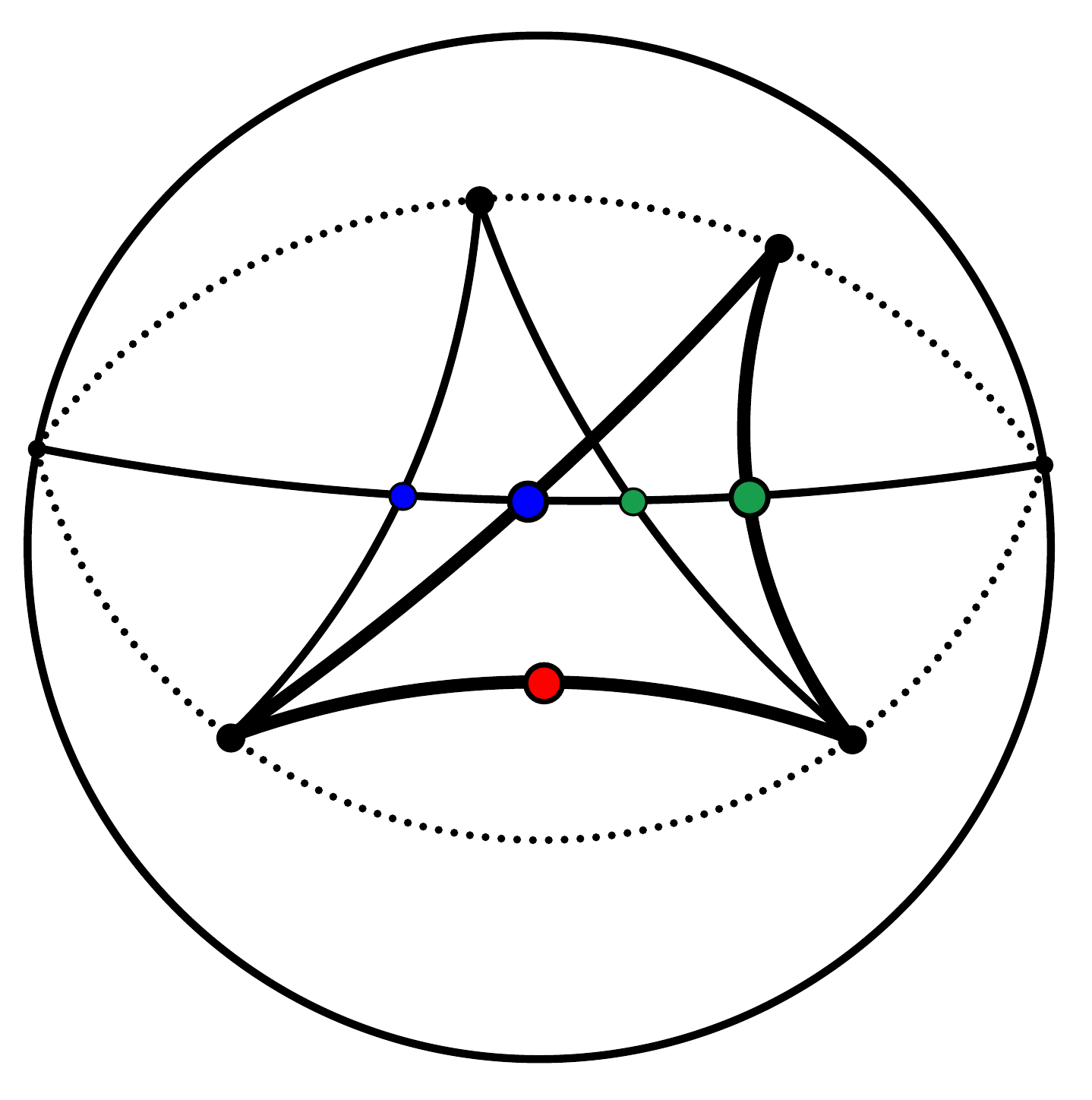}
\caption{Twist deformation: $X=Q\cdot R = Q'\cdot R'$.}
\label{f:earthquake}
\end{figure}

For this, let us fix $\theta\in(0,2\pi)$ and $\xi\in\fC_\theta$. We will use the notation from Figure~\ref{f:setup-small}
for the pseudo-fundamental domain $O_A O_C O_B R(O_C)$ of $\rho_\xi$.
Consider the equidistance locus $L$ to the geodesic $O_RO_Q$ containing the point $O_A$.
Take an arbitrary point $O_A'$ on $L$ and let \[\{O_R'\}=O_A'O_B\cap O_RO_Q,\quad \{O_Q'\}=O_A'O_C\cap O_RO_Q\] as in
Figure~\ref{f:earthquake}. Then we have the following.

\begin{lem}\label{l:twist}
Regardless of the position of $O_A'$ on $L$, we have that $d(O_R',O_Q')=d(O_R,O_Q)$,
and 
the areas of $O_AO_BO_C$ and $O_A'O_BO_C$ are the same.
\end{lem}

We will denote the rotation of angle $\alpha$ at the point $V\in\bH^2$ as $\Rot_\alpha(V)$. 

\bp[Proof of Lemma~\ref{l:twist}]
Let $Q'$ and $R'$ denote the $\pi$--rotations at $O_Q'$ and $O_R'$ respectively. 
The geodesic $O_RO_Q$ is the common invariant geodesic for the hyperbolic translations $Q\cdot R$ and $Q'\cdot R'$,
both of which map $O_B$ to $O_C$. Thus we have $Q\cdot R = Q'\cdot R'$ and hence $d(O_R',O_Q')=d(O_R,O_Q)$.
Moreover, we have that 
\[
\Rot_{\theta/2}(O_C)=C=P\cdot R\cdot Q=P\cdot R'\cdot Q'=\Rot_{\theta'/2}(O_C).\]
This implies that $\theta=\theta'$, which is the twice of the interior angle sum of each triangle.
\end{proof}

Recall our convention that a word $w\in W$ also denotes its image $\hat\rho_\xi(w)\in\PSL(2,\bR)$ when $\xi$ is fixed. The deformation from $\form{P,Q,R}$ to $\form{P,Q',R'}$ is called an \emph{$X$-twist}. Since the distance $h_X=d(R,Q)$ remains constant by the above lemma, we see from \eqref{e:trace} that the $X$-twist occurs along the coordinate curve $x=\textup{const.}=\tr X$.
Moreover, such a twist occurs inside $\fC_\theta$ since the interior angle sum is preserved.

Note that the distance $d_A=d(O_B,O_C)$ is invariant as well under $X$-twists.
In fact, the whole bi-infinite sequence of $\theta/2$-rotations $\{X^iBX^{-i}\}_{i\in\bZ}$
(as well as their various conjugates by $B$, $C$ and $P$) remain unchanged under $X$-twists. Lastly, observe that, by a suitable $X$-twist, we can always increase either $d(O_A,O_B)$ or $d(O_A,O_C)$ by an arbitrary amount.

Analogously, we define \emph {$Y$-twists} and \emph{$Z$-twists}.

\subsection{Torsion relations}
In this subsection we will show that holonomy maps containing torsion in their image are dense in $\fC_\theta$. For this, let us first recall  the following definition due to Lobachevskii.

Let $h\in\bR_{>0}$  be a positive  real number,  and let  $s\sse\bH^2$  be a geodesic segment of length $h$.  Choose a bi-infinite geodesic $\gamma$ which meets  an endpoint $p$ of $s$ at a right angle. We then connect the other  endpoint $q$ of  $s$ with an endpoint at  infinity of  $\gamma$ via a geodesic ray  $\delta$. The \emph{parallel angle} $\alpha(h)$ is defined as the acute angle at $q$ between $s$  and $\delta$. It is straightforward to  check  that $\alpha(h)$ is independent of all the choices that are  made.

\begin{thm}\label{t:torsion}
Let $\theta\in (0,2\pi)\setminus\bQ\pi$. Then there exists a dense subset $S_\theta\sse(2,\infty)$ such that 
every point $\xi\in \fC_\theta\cap (S_\theta\times\bR\times\bR)$ corresponds to the holonomy map $\rho_\xi$ satisfying
\[ 
\rho_\xi(r^m)=1\]
for some $r\in \pi_1(M)\setminus\{1\}$ and $m\ge2$. Analogous results also hold for
$\bR\times S_\theta\times \bR$ and $\bR\times\bR\times S_\theta$.
\end{thm}

\begin{figure}[ht]
\labellist
\pinlabel {\tiny $O_P$} at 210 200
\pinlabel {\tiny $O_B$} at 50 180
\pinlabel {\tiny $O_C$} at 370 180
\pinlabel {\tiny $O$} at 210 60
\endlabellist
\centering
\includegraphics[width=0.2\textwidth]{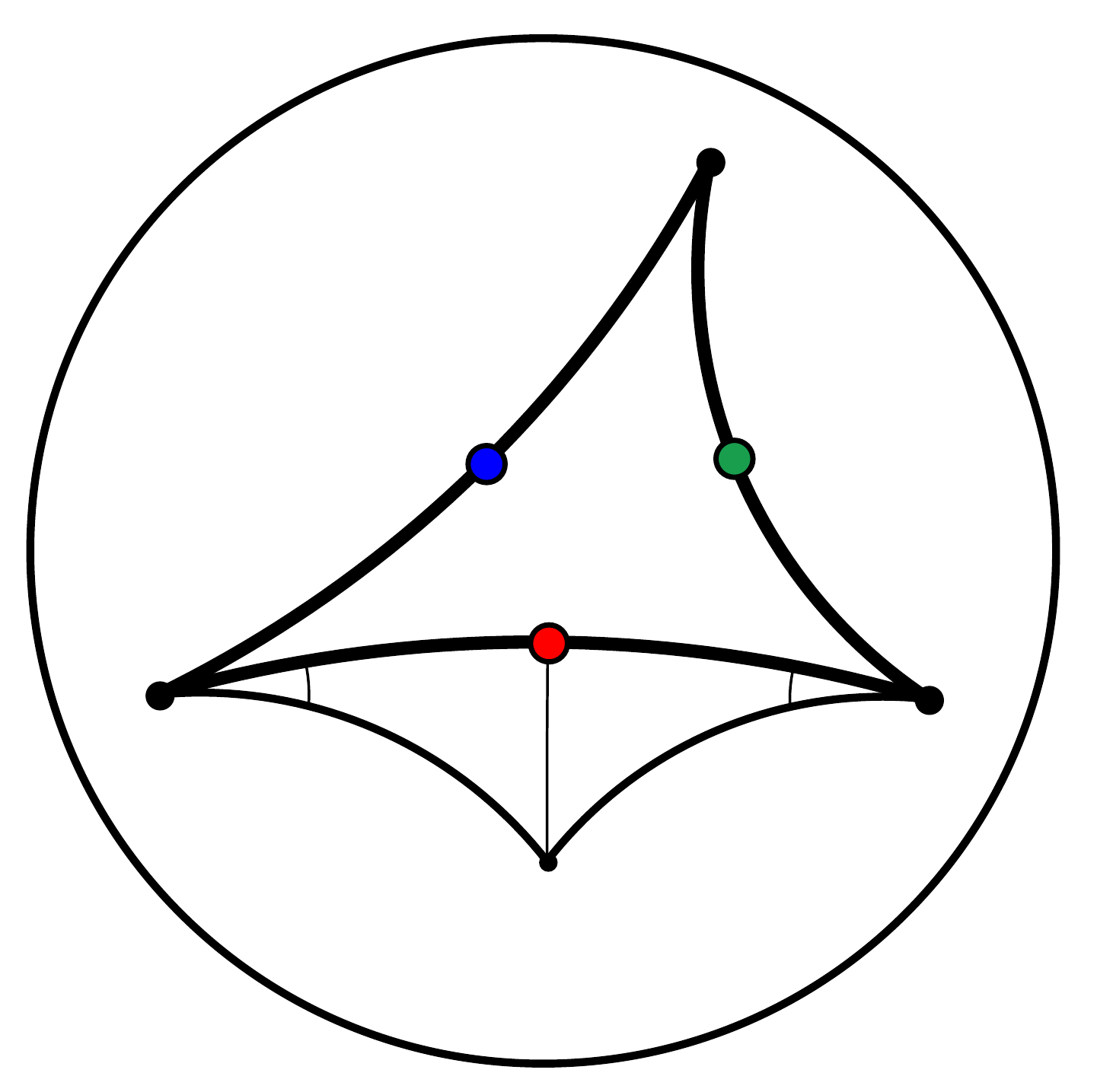}
\caption{Proof of Theorem~\ref{t:torsion}.}
\label{f:single}
\end{figure}

\bp[Proof of Theorem~\ref{t:torsion}]
We use the notation in Figure~\ref{f:setup-small}. Fix an arbitrary point $\xi\in\fC_\theta$.
Since $\theta\notin\bQ\pi$, we can find an integer $N>0$ so that the common rotation angle \[N\theta/2 \mod{2\pi}\] of $B^{N}$ and $C^{N}$ is
 between $0$ and $\alpha(d_A/2)$. Consider the isosceles triangle $OO_BO_C$ as in Figure~\ref{f:single} for which the
 common angle at $O_B$ and $O_C$ is $N\theta/2$.
 Let $J$, $J_B$, and $J_C$ denote the reflections along the geodesics $O_BO_C$, $OO_B$, and $OO_C$, respectively. Then we have
\begin{align*}
B^{2N}&=JJ_B,\\
C^{2N}&=J_CJ,\\
C^{2N}B^{2N}&=J_CJ_B=\Rot_\beta(O),
\end{align*}
where $-2\pi<\beta=-2\measuredangle(O_B,O,O_C)<0$. Since the common angle $N\theta/2$ at $O_B$ and $O_C$
is constant in $\fC_\theta$, the angle $\beta$ depends only on the distance $d_A$.

Now we perturb the representation (by a $Y$-twist or $Z$-twist, for example) in order to increase the distance $d_A$ slightly,
so that the deformed angle $\beta'$ satisfies $\beta'\in\bQ\pi$ and thus $C^{2N}B^{2N}=[Z,X]^N[Y,Z]^N$ becomes a torsion element.
Since $d_A$ is constant under $X$-twists, we have just found an element $s$ of $S_\theta$ such that every point in the $X$-coordinate curve
\[
\fC_\theta\cap(\{s\}\times\bR\times\bR)\]
yields a holonomy map with a torsion relation of the form $([Z,X]^N[Y,Z]^N)^m$ for some $m\ge2$.
Since we found such a curve which is arbitrarily close to any given point $\xi\in\fC_\theta$, the set $S_\theta$ must be dense in $(2,\infty)$.
\end{proof}

\subsection{Non-torsion-type kernel elements}
We now establish the following theorem, which obviously implies Theorem~\ref{t:nontorsion-intro}.

\begin{thm}\label{t:nontorsion}
If $\theta\in (0,2\pi)\setminus\bQ\pi$, then there exists a dense subset $S_\theta\sse(2,\infty)$ such that 
for each $s\in S_\theta$ every point $\xi$ on the coordinate curve 
\[\fC_\theta\cap \{z=s\}\]
corresponds to a holonomy map $\rho_\xi$ that has a palindromic non-torsion-type word in the kernel.
\end{thm}

By symmetry, the conclusion $z=s$ can be replaced by $x=s$ or $y=s$.

\begin{rem}\label{r:triple}
The theorem implies that there exists a dense subset 
\[ S_\theta' = \fC_\theta\cap ((\bR\times S_\theta\times S_\theta)\cup(S_\theta\times\bR\times S_\theta)\cup(S_\theta\times S_\theta\times\bR))\]
of $\fC_\theta$, consisting of points with two coordinates from $S_{\theta}$, such that each point $\xi\in S_\theta'$ corresponds to a holonomy map with non-torsion-type kernel elements. 
\end{rem}

\begin{figure}[ht]
\labellist
\pinlabel {\tiny $O_C$} at 190 475
\pinlabel {\tiny $O_A$} at 80 226
\pinlabel {\tiny $O_B$} at 245 220
\pinlabel {\tiny $O_R$} at 165 270
\pinlabel {\tiny $Z(O_R)=B^{-1}(O_R)$} at 490 290
\pinlabel {\tiny $B(O_R)$} at 390 200
\endlabellist
\centering
\includegraphics[width=0.3\textwidth]{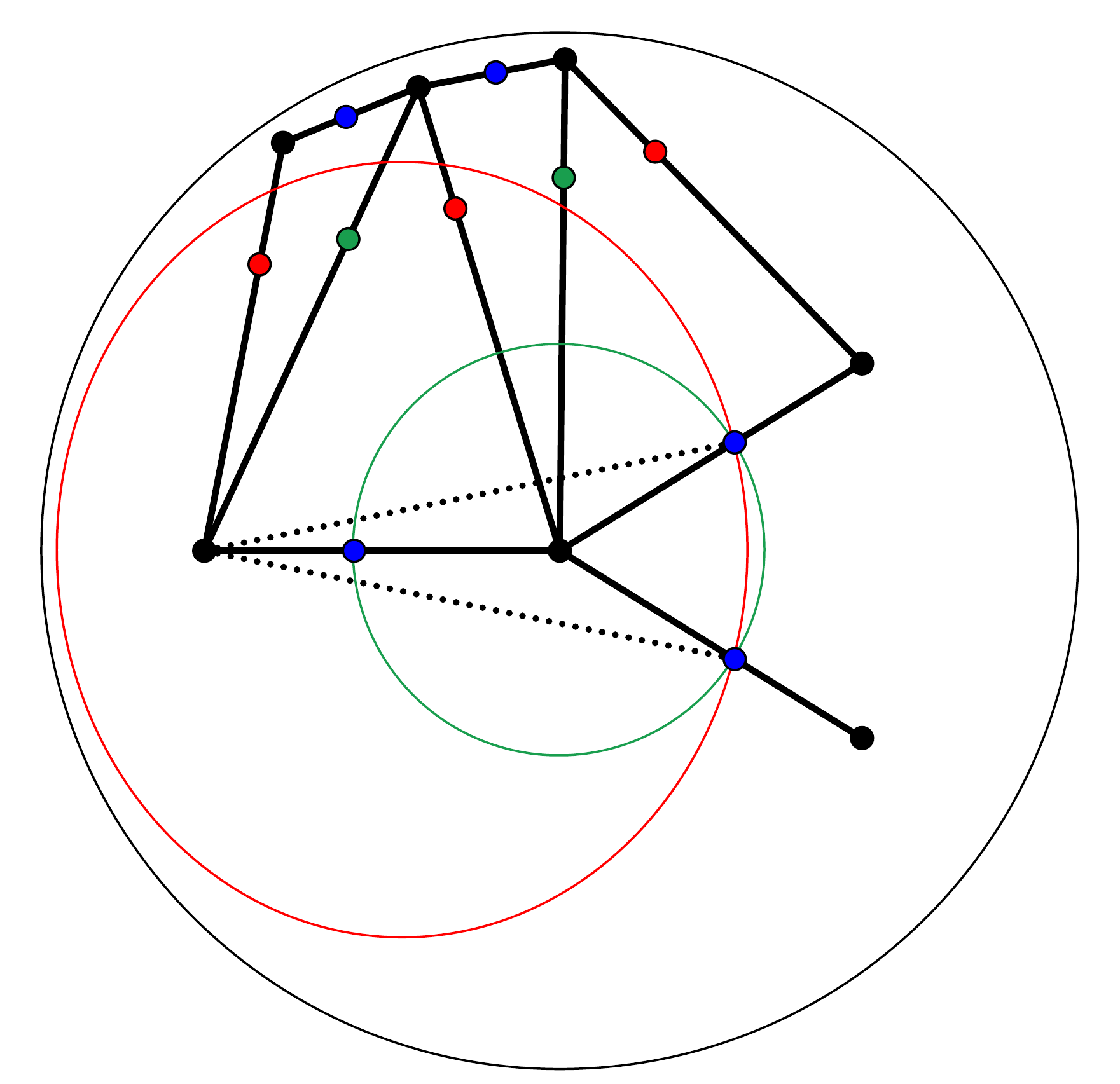}
\caption{Construction of the proof of Theorem~\ref{t:nontorsion}.}
\label{f:torsion_free_5}
\end{figure}

\bp[Proof of Theorem~\ref{t:nontorsion}]
We use the notation from Figure~\ref{f:setup-small} and from Section~\ref{s:coxeter}.
Since $Z(O_R)=B^{-1}(O_R)$, we have that
 \[d(O_A, Z(O_R))=d(O_A,B^{-1}(O_R))=d(O_A,B(O_R)).\]
In Figure~\ref{f:torsion_free_5}, two circles passing through $Z(O_R)$ and $B(O_R)$ are drawn
(using the Klein projective model of $\bH^2$); they are centered at $O_A$ and $O_B$, respectively. 

We claim that this configuration can be deformed by a suitable arbitrarily small $X$-twist so that 
two conjugates of the $\pi$-rotation $R$ coincide, that is,
\[
A^{2N}(BRB^{-1})A^{-2N}=ZRZ^{-1}.
\]

To see this claim, let us first recall that the angle \[\measuredangle(O_A,O_B,Z(O_R))=\measuredangle(O_A,O_B,B(O_R))=\theta/2\]
is constant under any deformation in $\fC_\theta$. Let \[\beta=\measuredangle(Z(O_R),O_A,B(O_R)),\] and we can choose $N>0$
so that the rotation angle \[N\theta=(2N)(\theta/2) \mod{2\pi}\] of $A^{2N}$ is arbitrarily close to but less than $\beta$.
Take a suitable  small $X$-twist
that increases the length $d(O_A,O_B)$; it will also increase the length \[d(O_B,O_R)=d(O_B,Z(O_R))=d(O_B,B(O_R)).\]
Then the angle $\beta$ decreases slightly and we twist until \[N\theta=\beta\mod{2\pi}.\]
(We note that this argument works for any $\theta\in(0,2\pi)\setminus {\mathbb Q}\pi$.) In the end, we obtain the relation
\[A^{2N}(BRB^{-1})A^{-2N}=ZRZ^{-1},\]
which is invariant under $Z$-twists. The proof of the claim is complete.

If we define
\[u_N:=[Z^{-1}A^{2N}B,R]\in F_2,\]
then $u_N$ is a product of two involutions and the above relation is equivalent to  $\rho_\xi(u_N)=1$.
Thus for a dense choice of $s$ in $(2,\infty)$ there exists an integer $N=N(s)>0$ such that
the holonomy map of every point $\xi$ on the coordinate curve $\{z=s\}$ satisfies  $u_N\in\ker\rho_\xi$. 

We can write $u_N$ as a word in $F_2$ from the computation below.
\begin{align*}
u_N
=&(Z^{-1}A^{2N}B)R(B^{-1}A^{-2N}Z)R\\
=&Z^{-1}[X,Y]^N(RQP)R(PQR)(RPQRPQ)^N(PQ)R\\
=&Z^{-1}[X,Y]^N (RQ) (PR) (PQ) [R(RPQRPQ)^N P]QR\\
=&Z^{-1}[X,Y]^N X^{-1} Y^{-1}Z[R(R)(PQRPQR)^{N-1}(PQRPQ)P]X\\
=&Z^{-1}[X,Y]^NX^{-1}Y^{-1}Z(ZYX)^{N-1}ZYZ^{-1}X\\
=&XY[X,Y]^N  X^{-1}Y^{-1}Y^{-1}X^{-1}[Y^{-1},X^{-1}]^{N-1}Y^{-1}X^{-1}YXYX\\
=&XY[X,Y]^N  X^{-1}Y^{-2}X^{-1}[Y^{-1},X^{-1}]^{N}YX.
\end{align*}

As can be expected from Lemma~\ref{l:pal}, $u_N$ is indeed palindromic (but not a proper power).
To finish the proof it suffices to have the following proposition, which is proved in the next section.
\ep

\begin{prop}\label{p:palindrome}
For each $N\ge1$, we have that $u_N$ is of non-torsion-type.
\end{prop}

\begin{rem}\label{rem:palindrome}
In the above proof, by Proposition~\ref{p:equiv}, the image $\rho_\xi(u_N)$ is trivial for $\xi\in\fC_\theta\cap\{z=s\}$,
while it is a nontrivial hyperbolic translation for  $\xi\in\fC_\theta\setminus\{z=s\}$.
\end{rem}

Recall we denote by $g_w(x,y,z)\in\bZ[x,y,z]$ the trace polynomial of a word $w\in F_2$. From the above proof and the remark following it, we see for each $s\in S_\theta$ there exists a positive integer $N=N(\theta,s)$ such that the coordinate curve $\fC_\theta\cap\{z=s\}$ can be expressed as the intersection of two surfaces\[\fC_\theta\cap\{z=s\}=\fC_\theta\cap\{g_{u_N}=\pm2\}.\]

In general, two words $U$ and $V$ may admit the same trace polynomial. For instance, we have $g_U=g_V$ if $U$ and $V$ are conjugate, or if $U = s_1 \ldots s_\ell$ and $V=s_\ell s_{\ell-1}\cdots s_1$ for $s_i\in\{X^{\pm1},Y^{\pm1}\}$. In our case, the word $u_N$ is palindromic and so, it seems rare that another surface $\{g_w=\pm2\}$ contains the above coordinate curve for a different word $w\in F_2$ not conjugate to $u_N$. This raises the following question:

\begin{que} Let $\theta\in (0,2\pi)\setminus\bQ\pi$. Does there exist a dense subset $R \subset S_\theta$ such that for each $r\in R$ and for $N=N(\theta,r)$ as above, a very general point $\xi$ in the coordinate curve $\fC_\theta\cap\{z=r\}$ corresponds to the monodromy $\rho_\xi$ whose image is isomorphic to the one-relator group \[\form{X,Y\mid u_N=1}?\] \end{que}

\section{Proof of Proposition~\ref{p:palindrome}}

In this section, we prove Proposition~\ref{p:palindrome}. The proof  we offer
 is somewhat technical by nature, but based on a beautiful solution by Newman of the word problem 
 in one-relator groups with torsion (Corollary~\ref{c:dehn}) as we describe now.

Magnus (1932) discovered that the word problem for a (general) one-relator group is solvable.
In the case of a one-relator with torsion, the following theorem due to Newman (1968) gives a particularly simple solution.
For a reduced word $w$ in a free group, we let $|w|$ denote its word length and we let $[w]$ denote its cyclic conjugacy class.
\begin{thm}[Newman Spelling Theorem]\label{t:newman}
Let $G = \form{x_1,\ldots,x_n\mid r^k}$ be a one-relator group for some cyclically reduced word $r$ in $F_n:=F(x_1,\ldots,x_n)$,
where  $k\ge2$. 
Suppose that a nontrivial reduced word $w\in F_n$ belongs to the kernel of the natural quotient map $F_n\to G$.
Then $w$ contains a subword $u$, of length strictly larger than $(k-1)|r|$, such that $u$ is a subword of $r^{\pm k}$.
\end{thm}

 If $u$ is a subword of $v$, we write $u\preccurlyeq v$.
By the Newman Spelling Theorem, we have the following solution to the word problem in two-generator one-relator group with torsion, which
 we call the \emph{Dehn--Newman algorithm}, as it is based on the Dehn's solution to the word problem for surface groups. 
Recall our convention that $F_2=\form{X,Y}$.

\begin{cor}[Dehn--Newmann Algorithm]\label{c:dehn}
Let $u$ and $r$ be cyclically reduced words in $F_2$, and let $m\ge2$.
Then the truth value of the statement $u\in\fform{r^m}$ can be decided in the following steps, within a finite time.
\be
\item[\textbf {Step 1}]\label{p:dehn1} If there exists a cyclic conjugation $R=R_1R_2$ of $r^{\pm1}$ such that $|R_1|=1$
and such that  $R^{m-1}R_1$ is a subword of $u$, then proceed to \textbf{Step 2}; otherwise, conclude that $u\not\in\fform{r^m}$.
\item[\textbf {Step 2}]\label{p:dehn2}  Let $U$ be the word obtained from $u$ by substituting $R_2^{-1}$ for $R^{m-1}R_1$.
If $U\ne1$, then repeat \textbf{Step 1}; otherwise, conclude that $u\in\fform{r^m}$.
\ee
\end{cor}

Let us first see a simple application of the algorithm.
\begin{lem}\label{l:un1}
Let $N\ge1$. 
\be
\item\label{p:un1-1}
If $r\in \{Y,XY,[X,Y]\}$ and $m\ge2$, then $u_N\not\in\fform{r^m}$.
\item\label{p:un1-2}
If $r\in F_2$ and $m\ge3$, then $u_N\not\in\fform{r^m}$.
\ee
\end{lem}
\bp
(\ref{p:un1-1})
We let $G:=F_2/\fform{r^m}$. Assume for contradiction that $u_N=1$ in $G$.

Suppose first that $r=Y$. By
the Dehn--Newman Algorithm (Corollary~\ref{c:dehn}), we should have $Y^{\pm m}\preccurlyeq u_N$.
This implies $m=2$. In $G=F_2/\fform{Y^2}$, the element $u_N$ coincides with the word
\[
u':=XY[X,Y]^NX^{-2}[Y^{-1},X^{-1}]^NYX.\]
Since $Y^{\pm2}\not\preccurlyeq u'$, we see that $u_N=u'\ne1$ in $G$. 

Now, consider the case $r=XY$. Since a cyclic conjugation of a length five subword of $(XY)^3$ does not appear
in $u_N$ as a subword, we see that  $m=2$. We cancel out all the occurrences of 
\[
(XY)^{\pm 2}, (YX)^{\pm2}\]
in $u_N$ to obtain 
\[
u':=X^{-1}Y^{-1}[X,Y]^{N-2}XYYX[Y^{-1},X^{-1}]^{N-2}Y^{-1}X^{-1}.\]
Since $u'$ does not contain a length three subword of 
\[
(XY)^{\pm 2}, (YX)^{\pm2}\]
we conclude from the Dehn--Newman algorithm that $u'\ne1$ in $G$.

Finally, assume $r=[X,Y]$. Writing $N=mq+t$ for $t\in[0,m-1]$, the element $u_N$ can be rewritten in $G$ as
\[
u':=XY[X,Y]^t X^{-1}Y^{-2}X^{-1}[Y^{-1},X^{-1}]^tYX.\]
It is obvious that $u'$ does not contain a subword of $r^m$ with length $4(m-1)+1$. Hence, we see $u'\ne1$ in $G$.

(\ref{p:un1-2})
Suppose not. The Dehn--Newman Algorithm implies that
for some cyclically reduced word $r$ that is not a proper power 
we have $r^2\preccurlyeq u$.
This consideration restricts the choices of $r$ to be a cyclic conjugation of a word in 
\[\{Y, XY,[X,Y]\}^{\pm1}.\]
By part (\ref{p:un1-1}), this is impossible.
\ep

We are now ready to give a proof of the proposition.

\bp[Proof of Proposition~\ref{p:palindrome}]
Suppose  $u_N\in\fform{r^m}$
for some $r\in F_2$ and $m\ge2$.
By Lemma~\ref{l:un1}, we need only consider the case $m=2$.
Consider the  following procedure, which takes $u=u_N$  as  an input.
\be[(i)]
\item Pick a cyclically reduced word $r=r_1 r_2$ in $F_2$ such that $|r_1|=1$ and such that $r_1 r_2 r_1\preccurlyeq u$;
here, we require that $r_1r_2$ is not cyclically conjugate to a word in 
\[\{Y,XY,[X,Y]\}^{\pm1}.\]
Proceed to Step (ii).
\item Replace an occurrence of $r_1 r_2  r_1$ in $u$ by $r_2^{-1}$ to obtain a reduced word $u'$. Proceed to Step (iii).
\item Pick a cyclic conjugation $r_1' r_2'$ of $r^{\pm1}$ such that $|r_1'|=1$ and such that $r_1' r_2' r_1'\preccurlyeq u'$.
Plug  $u\leftarrow u'$ and $r_i\leftarrow r_i'$. Repeat Step (ii).
\ee
This process is a simplified version of Dehn--Newman Algorithm. It is obvious that if  the procedure
 terminates at step (i) or (iii), we can conclude that $u_N$ is of non-torsion-type. 

We can enumerate all the possible choices of $r_1r_2r_1$ in step (i) as follows, excluding the obvious case $r_1r_2r_1=u_N$.
Verifying all of the following cases,  we can conclude that $u_N\not\in\fform{r^m}$ for all $r\in F_2$ and $m\ge2$.
\be
\item\label{p:X1} $XY[X,Y]^tX$ for $1\le t\le N-1$;
\item\label{p:X2}  $XY[X,Y]^NX^{-1}Y^{-2}X^{-1}[Y^{-1},X^{-1}]^t$ for $1\le t\le N$;
\item\label{p:X3} $[X,Y]^sX^{-1}Y^{-2}X^{-1}[Y^{-1},X^{-1}]^t$ for $1\le s,t\le N$;
\item\label{p:X4}  $[X,Y]^sX^{-1}Y^{-2}X^{-1}[Y^{-1},X^{-1}]^NYX$ for $1\le s\le N$;
\item\label{p:X5}  $X[Y^{-1},X^{-1}]^tYX$ for $1\le t\le N-1$;
\item\label{p:X11}  $X^{-1}Y^{-1}[X,Y]^tX^{-1}$ for $1\le t\le N-1$;
\item\label{p:X12}  $X^{-1}Y^{-1}[X,Y]^tX^{-1}Y^{-2}X^{-1}$ for $0\le t\le N-1$;
\item\label{p:X13}  $X^{-1}Y^{-1}[X,Y]^sX^{-1}Y^{-2}X^{-1}[Y^{-1},X^{-1}]^t Y^{-1}X^{-1}$ for $0\le s, t\le N-1$;
\item\label{p:X14}  $X^{-1}Y^{-2}X^{-1}$;
\item\label{p:X15}  $X^{-1}Y^{-2}X^{-1}[Y^{-1},X^{-1}]^t Y^{-1}X^{-1}$ for $0\le t\le N-1$;
\item\label{p:X16} $X^{-1}[Y^{-1},X^{-1}]^t Y^{-1}X^{-1}$ for $1\le t\le N-1$;
\item\label{p:Y1} $Y[X,Y]^tXY$ for $1\le t\le N-1$;
\item\label{p:Y2} $Y[X,Y]^NX^{-1}Y^{-2}X^{-1}[Y^{-1},X^{-1}]^tY^{-1}X^{-1}Y$ for $0\le t\le N-1$;
\item\label{p:Y3} $Y[X,Y]^NX^{-1}Y^{-2}X^{-1}[Y^{-1},X^{-1}]^NY$;
\item\label{p:Y4} $YX^{-1}Y^{-1}[X,Y]^sX^{-1}Y^{-2}X^{-1}[Y^{-1},X^{-1}]^tY^{-1}X^{-1}Y$ for $0\le s, t\le N-1$;
\item\label{p:Y5} $YX^{-1}Y^{-1}[X,Y]^sX^{-1}Y^{-2}X^{-1}[Y^{-1},X^{-1}]^NY$ for $0\le s, t\le N-1$;
\item\label{p:Y6} $YX[Y^{-1},X^{-1}]^tY$ for $1\le t\le N-1$;
\item\label{p:Y11} $Y^{-1}[X,Y]^tX^{-1}Y^{-1}$ for $1\le t\le N-1$;
\item\label{p:Y12} $Y^{-1}[X,Y]^tX^{-1}Y^{-2}$ for $0\le t\le N-1$;
\item\label{p:Y13} $Y^{-1}[X,Y]^sX^{-1}Y^{-2}X^{-1}[Y^{-1},X^{-1}]^t Y^{-1}$ for $0\le s, t\le N-1$;
\item\label{p:Y14} $Y^{-2}X^{-1}[Y^{-1},X^{-1}]^t Y^{-1}$ for $0\le t\le N-1$;
\item\label{p:Y15} $Y^{-1}X^{-1}[Y^{-1},X^{-1}]^t Y^{-1}$ for $1\le t\le N-1$.
\ee
The following pairs of cases yield the same relator $r^2$:
\begin{itemize}
\item (\ref{p:X1}) and (\ref{p:Y1});
\item (\ref{p:X11}) and (\ref{p:Y11});
\item (\ref{p:X16}) and (\ref{p:Y15}).
\end{itemize}

So, for each of the 19 cases remaining, we verify that the above algorithm terminates.

{\bf Case (\ref{p:X1}).} This case corresponds to 
\[
r_1=X,\quad r_2=Y[X,Y]^t\text{ for }1\le t\le N.\]
Replacing $r_1r_2r_1=XY[X,Y]^tX$ in $u_N$ by $r_2^{-1}=[Y,X]^tY^{-1}$, we obtain
\begin{align*}
u'&= [Y,X]^t X^{-1}Y^{-1} [X,Y]^{N-t-1}X^{-1}Y^{-2}X^{-1}[Y^{-1},X^{-1}]^{N}YX
\\
&=[Y,X]^t X^{-1}Y^{-1} [X,Y]^{N-t-1}X^{-1}Y^{-2}
X^{-1}
(Y^{-1}X^{-1})
[Y,X]^{t}
[Y,X]^{N-t-1}(YX)^2
\end{align*}
We now continue with 
\[
r_1'=X^{-1},\quad r_2'=Y^{-1}X^{-1}[Y,X]^{t-1}YXY^{-1}.\]
One sees that
\[
r_1'r_2'r_1'=X^{-1}Y^{-1}X^{-1}[Y,X]^t\preccurlyeq u'.\]
Plug $u\leftarrow u'$ and $r_i\leftarrow r_i'$. We replace $r_1r_2r_1$ in $u$ by $r_2^{-1}$
to obtain
\begin{align*}
u'&:=[Y,X]^t X^{-1}Y^{-1} [X,Y]^{N-t-1}X^{-1}Y^{-1}
X^{-1}Y^{-1}[X,Y]^{t-1}XY
[Y,X]^{N-t-1}(YX)^2
\end{align*}
It is then clear that Step (iii) cannot be executed for $u'$ and $r$.
We conclude that $u'=u\ne1$ in $F_2/\fform{r^2}$.

{\bf Case (\ref{p:X2}).}
In this case we have 
\[
r_1=X,\quad r_2=Y[X,Y]^NX^{-1}Y^{-2}X^{-1}[Y^{-1},X^{-1}]^{t-1}Y^{-1}X^{-1}Y\text{ for }1\le t\le N.\]
Then $r_2^{-1}=Y^{-1}XY[X^{-1},Y^{-1}]^{t-1}XY^2X[Y,X]^NY^{-1}$. We replace $r_1r_2r_1$ in $u_N$ by  $r_2^{-1}$ to obtain
\[u'=Y^{-1}XY[X^{-1},Y^{-1}]^{t-1}XY^2X[Y,X]^NY^{-1}[Y^{-1},X^{-1}]^{N-t}YX,\] which is reduced.
It is clear that no cyclic conjugate of $(r_1r_2)^{\pm 1}$ occurs as
a subword of  $u'$, as is seen from aligning $[Y,X]^N$ in $u'$ with the corresponding copy in $r_2^{-1}$.

{\bf Case (\ref{p:X3}).}
In this case we have 
\[
r_1=X,\quad r_2=YX^{-1}Y^{-1}[X,Y]^{s-1}X^{-1}Y^{-2}X^{-1}[Y^{-1},X^{-1}]^{t-1}Y^{-1}X^{-1}Y\text{ for }1\le s,t\le N.\] Then
$r_2^{-1}=Y^{-1}XY[X^{-1},Y^{-1}]^{t-1}XY^2X[Y,X]^{s-1}YXY^{-1}$. Replacing $r_1r_2r_1$ by $r_2^{-1}$ gives us
\[u'=XY[X,Y]^{N-s}Y^{-1}XY[X^{-1},Y^{-1}]^{t-1}XY^2X[Y,X]^{s-1}YXY^{-1}[Y^{-1},X^{-1}]^{N-t}YX.\] This word is reduced unless
$N=s$ or $N=t$, in which  case a subword $X^2$ appears.
It is clear that no 
cyclic conjugate of $(r_1r_2)^{\pm 1}$ occurs as
a subword of  $u'$.

{\bf Case (\ref{p:X4}).}
In this case we have 
\[
r_1=X,\quad r_2=YX^{-1}Y^{-1}[X,Y]^{s-1}X^{-1}Y^{-2}X^{-1}[Y^{-1},X^{-1}]^{N}Y\text{ for }1\le s\le N.\] Then
$r_2^{-1}=Y^{-1}[X^{-1},Y^{-1}]^NXY^2X[Y,X]^{s-1}YXY^{-1}$. Replacing $r_1r_2r_1$ by $r_2^{-1}$ gives us
\[
u'=XY[X,Y]^{N-s}Y^{-1}[X^{-1},Y^{-1}]^NXY^2X[Y,X]^{s-1}YXY^{-1},
\]
which is reduced. There is exactly one occurrence of $Y^2$ in $u'$ and in $(r_1r_2)^{-1}$, and if these are aligned then a further
reduction can be made. Clearly such an alignment cannot be made, as is checked by reading to the left or right from $Y^2$ in $u'$.
Thus, we must be able to find a subword of $u'$ which is a cyclic conjugate of $(r_1r_2)^{\pm1}$ and which begins and ends with the
same power of $Y$, that is \[Y^{-1}X^{-1}[Y^{-1},X^{-1}]^{N}YXYX^{-1}Y^{-1}[X,Y]^{s-1}X^{-1}Y^{-1}\] or its inverse. It is clear that this
does not happen.

{\bf Case (\ref{p:X5}).}
In this case  we  have
\[
r_1=X,\quad r_2=[Y^{-1},X^{-1}]^tY \text{ for }1\le t\le N-1.\] 
Then $r_2^{-1}=Y^{-1}[X^{-1},Y^{-1}]^t$.
Replacing $r_1r_2r_1$ by $r_2^{-1}$ gives us 
\[
u'=XY[X,Y]^t X^{-1}Y^{-2} X^{-1}  [Y^{-1},X^{-1}] Y^{-1}X^{-1}YY^{-1}[X^{-1}, Y^{-1}]^t,
\]

which reduces  to \[XY[X,Y]^tX^{-1}Y^{-2}X^{-1}[Y^{-1},X^{-1}]Y^{-1}X^{-1}[X^{-1}, Y^{-1}]^t.\]

This latter  word is both
reduced and cyclically reduced. In order to perform further  reductions, $u'$ must contain subwords of the form
\[YX[Y^{-1},X^{-1}]^t\quad\textrm{or}\quad X^{-1}Y^{-1}[X^{-1},Y^{-1}]^t.\]
For the first of  these  cases to occur  we would need $t=1$, and  then it is  clear  that
no cyclic  permutation of $u'$ contains  a cyclic permutation of this word.
 In the second  case, we similarly see that  no  cyclic conjugate of this word  occurs as a subword of a cyclic permutation of  $u'$. 

{\bf Case (\ref{p:X11}).}
In this case we have 
\[
r_1=X^{-1},\quad r_2=Y^{-1}[X,Y]^t \text{ for }1\le t\le N-1.\] 
Then
$r_2^{-1}=[Y,X]^tY$. Replacing $r_1r_2r_1$ by $r_2^{-1}$ gives us
\[
u'=XY[X,Y]^{N-t-1}XY[Y,X]^tY^3X^{-1}[Y^{-1},X^{-1}]^NYX,
\]
which is reduced. In order  to perform further reductions, we must align $[X,Y]^tX^{-1}Y^{-1}$ or
$YX[Y,X]^t$ in $u'$ and $(r_1r_2)^{\pm}$.
By inspection, this is not possible.

{\bf Case (\ref{p:X12}).}
In this case we have 
\[
r_1=X^{-1},\quad r_2=Y^{-1}[X,Y]^tX^{-1}Y^{-2} \text{ for }0\le t\le N-1.\] Then
$r_2^{-1}=Y^{2}X[Y,X]^tY$. Replacing $r_1r_2r_1$ by $r_2^{-1}$ gives us
\[
u'=XY[X,Y]^{N-t-1}XY^{3}X[Y,X]^tY[Y^{-1},X^{-1}]^NYX,
\]
which reduces to \[XY[X,Y]^{N-t-1}XY^{3}X[Y,X]^tX^{-1}YX[Y^{-1},X^{-1}]^{N-1}YX.\] 
There is only one power of $Y$ of absolute value more than one in $u'$, and it is clear that \[Y^{2}X[Y,X]^tYX=(r_1r_2)^{-1}\] does not
occur as a subword of $u'$ or of its cyclic conjugates, nor do any of the cyclic conjugates of $(r_1r_2)^{-1}$ which keep the
$Y^2$ intact. We therefore consider subwords
of the form $YX[Y,X]^tYXY$ and its inverse. It is immediate to verify that there are no such subwords, even in the degenerate case
where $t=0$.

{\bf Case (\ref{p:X13}).}
In this case we have 
\[
r_1=X^{-1},\quad r_2=Y^{-1}[X,Y]^sX^{-1}Y^{-2}X^{-1}[Y^{-1},X^{-1}]^t Y^{-1} \text{ for }0\le s,t\le N-1.\] Then
$r_2^{-1}=Y[X^{-1},Y^{-1}]^tXY^2X[Y,X]^sY$. Replacing $r_1r_2r_1$ by $r_2^{-1}$ gives us
\[
u'=XY[X,Y]^{N-t-1}XY^2[X^{-1},Y^{-1}]^tXY^2X[Y,X]^sY^2X[Y^{-1},X^{-1}]^{N-t-1}YX,
\]
which is reduced and cyclically reduced. Up to cyclic permutation, we must align
subwords of $(r_1r_2)^{\pm 1}$ and $u'$ of the form $X^{-1}Y^{-2}X^{-1}$
or $XY^2X$, or we must set $r_1'=Y$ and \[r_1'r_2'r_1'=YX[Y,X]^sYXY[X^{-1},Y^{-1}]^tXY^2.\]
Reading
to the left or to the right in $u'$ from this copy of $Y^2$ quickly yields a contradiction.

{\bf Case (\ref{p:X14}).}
In this case we have 
\[
r_1=X^{-1},\quad r_2=Y^{-2}.\] Then
$r_2^{-1}=Y^2$. Replacing $r_1r_2r_1$ by $r_2^{-1}$ gives us
\[
u'=XY[X,Y]^NY^2[Y^{-1},X^{-1}]^NYX=XY[X,Y]^{N-1}XYX^{-2}YX[Y^{-1},X^{-1}]^{N-1}YX,
\]
where the latter expression is both reduced and cyclically reduced. In order to reduce further, the word $u'$ must have a subword of
the form $Y^2$ or $Y^{-2}$. Indeed, such a subword would arise from either $(r_1r_2)^{\pm1}$ or $(r_2r_1)^{\pm1}$ as a subword of
$u'$, or by setting $r_1'=Y^{\pm1}$ and \[r_1'r_2'r_1'\in\{YXY^2,Y^{-1}X^{-1}Y^{-2}\}.\] Clearly none of these cases occur.

{\bf Case (\ref{p:X15}).}
In this case we have 
\[
r_1=X^{-1},\quad r_2=Y^{-2}X^{-1}[Y^{-1},X^{-1}]^tY^{-1}\text{ for }0\le t\le N-1.\] Then
$r_2^{-1}=Y[X^{-1},Y^{-1}]^tXY^2$. Replacing $r_1r_2r_1$ by $r_2^{-1}$ gives us
\[
u'=XY[X,Y]^NY[X^{-1},Y^{-1}]^tXY^2YX[Y^{-1},X^{-1}]^{N-t-1}YX,
\]
which reduces  to \[XY[X,Y]^{N-1}XYX^{-1}[X^{-1},Y^{-1}]^tXY^3X[Y^{-1},X^{-1}]^{N-t-1}YX,\] which is  both reduced and cyclically reduced.
Up to cyclic permutation, we must have that $u'$ has a subword of the form $X^{-1}Y^{-2}X^{-1}$ or $XY^2X$, or a cyclic permutation
of a word of the form
\[r_1'r_2'r_1'=YXY[X^{-1},Y^{-1}]^tXY^2\quad\textrm{or}\quad r_1'r_2'r_1'=Y^{-1}X^{-1}[Y^{-1}, X^{-1}]Y^{-1} X^{-1} Y^{-2}\] 
in order to perform further
reductions. Clearly this is not the case.

{\bf Case (\ref{p:X16}).}
In this case we have 
\[
r_1=X^{-1},\quad r_2=[Y^{-1},X^{-1}]^tY^{-1}\text{ for }1\le t\le N-1.\] 
Then
$r_2^{-1}=Y[X^{-1},Y^{-1}]^t$. Replacing $r_1r_2r_1$ by $r_2^{-1}$ gives us
\[
u'=XY[X,Y]^NX^{-1}Y^{-2}Y[X^{-1},Y^{-1}]^tYX[Y^{-1},X^{-1}]^{N-t-1}YX,
\]
which reduces to \[XY[X,Y]^NX^{-1}Y^{-1}[X^{-1},Y^{-1}]^tYX[Y^{-1},X^{-1}]^{N-t-1}YX.\] This last expression is  both reduced and
cyclically reduced. In  order to perform further reductions, a cyclic permutation of $u'$ must contain a subword which is
cyclic permutation of a word of the form
$XY[X^{-1},Y^{-1}]^t$ or $[Y^{-1},X^{-1}]^tY^{-1}X^{-1}$. The first of these does not, since $t>0$.

Thus we are reduced to the case where
 \[
r_1'=Y^{-1},\quad r_2'=X^{-1}YX[Y^{-1},X^{-1}]^{t-1}Y^{-1}X^{-1}.\]
Then, \[r_1'r_2'r_1'=[Y^{-1},X^{-1}]^tY^{-1}X^{-1}Y^{-1},\] which does not occur as a subword of $u'$ up to cyclic permutation.

{\bf Case (\ref{p:Y2}).}
In this case we have 
\[
r_1=Y,\quad r_2=[X,Y]^NX^{-1}Y^{-2}X^{-1}[Y^{-1},X^{-1}]^tY^{-1}X^{-1}\text{ for }0\le t\le N-1.\] Then
$r_2^{-1}=XY[X^{-1},Y^{-1}]^tXY^2X[Y,X]^N$. Replacing $r_1r_2r_1$ by $r_2^{-1}$ gives us
\[
u'=X^2Y[X^{-1},Y^{-1}]^tXY^2X[Y,X]^NX[Y^{-1},X^{-1}]^{N-t-1}YX,
\]
which reduces to \[X^2Y[X^{-1},Y^{-1}]^tXY^2X[Y,X]^{N-1}YXY^{-1}[Y^{-1},X^{-1}]^{N-t-1}YX.\] This last expression is reduced and
cyclically reduced. In order to perform further reductions, we must have that after performing a cyclic permutation,
there is a subword of $u'$ of the form  \[Y[X,Y]^NX^{-1}Y^{-2}X^{-1}[Y^{-1},X^{-1}]^tY^{-1}X^{-1}=r_1r_2\] or its inverse. By considering
the necessary alignment of product of commutators $[Y,X]^N$, we quickly obtain a contradiction.

{\bf Case (\ref{p:Y3}).}
In this case we have 
\[
r_1=Y,\quad r_2=[X,Y]^NX^{-1}Y^{-2}X^{-1}[Y^{-1},X^{-1}]^N.\] Then
$r_2^{-1}=[X^{-1},Y^{-1}]^NXY^2X[Y,X]^N$. Replacing $r_1r_2r_1$ by $r_2^{-1}$ gives us
\[
u'=X[X^{-1},Y^{-1}]^NXY^2X[Y,X]^NX,
\]
which reduces to \[Y^{-1}XY[X^{-1},Y^{-1}]^{N-1}XY^2X[Y,X]^NX,\] which is both reduced and cyclically reduced.  Note that the words
$(r_1r_2)^{\pm}$ have length $4N+5$ and that $u'$ has length $4N+4$,
whence no cyclic  permutation of $(r_1r_2)^{\pm}$ occurs as a subword of  $u'$.

{\bf Case (\ref{p:Y4}).}
In this case we have 
\[
r_1=Y,\quad r_2=X^{-1}Y^{-1}[X,Y]^sX^{-1}Y^{-2}X^{-1}[Y^{-1},X^{-1}]^tY^{-1}X^{-1}\text{ for }0\le s,t\le N-1.\] Then
$r_2^{-1}=XY[X^{-1},Y^{-1}]^tXY^2X[Y,X]^sYX$. Replacing $r_1r_2r_1$ by $r_2^{-1}$ gives us
\[
u'=XY[X,Y]^{N-s-1}X^2Y[X^{-1},Y^{-1}]^tXY^2X[Y,X]^sYX^2[Y^{-1},X^{-1}]^{N-t-1}YX,
\]
which is both reduced and cyclically reduced. In  order to perform further reductions, we would require a cyclic permutation of 
\[Y^{-1}XY[X^{-1},Y^{-1}]^tXY^2X[Y,X]^sYX\] or its inverse to occur as a subword of  $u'$, which it does not.

{\bf Case (\ref{p:Y5}).}
In this case we have 
\[
r_1=Y,\quad r_2=X^{-1}Y^{-1}[X,Y]^sX^{-1}Y^{-2}X^{-1}[Y^{-1},X^{-1}]^N\text{ for }0\le s\le N-1.\] Then
$r_2^{-1}=[X^{-1},Y^{-1}]^NXY^2X[Y,X]^sYX$. Replacing $r_1r_2r_1$ by $r_2^{-1}$ gives us
\[
u'=XY[X,Y]^{N-s-1}X[X^{-1},Y^{-1}]^NXY^2X[Y,X]^sYX^2,
\]
which reduces to 
\[
XY[X,Y]^{N-s-1}Y^{-1}XY[X^{-1},Y^{-1}]^{N-1}XY^2X[Y,X]^sYX^2.
\]
This last expression is reduced and cyclically reduced. In order to perform further reductions, we must have a subword of $u'$ which
is a cyclic permutation of \[Y^{-1}[X^{-1},Y^{-1}]^NXY^2X[Y,X]^sYX\] or its inverse. This is  not the  case.

{\bf Case (\ref{p:Y6}).}
In this case we have 
\[
r_1=Y,\quad r_2=X[Y^{-1},X^{-1}]^t\text{ for }0\le t\le N-1.\] Then
$r_2^{-1}=[X^{-1},Y^{-1}]^tX^{-1}$. Replacing $r_1r_2r_1$ by $r_2^{-1}$ gives us
\[
u'=XY[X,Y]^NX^{-1}Y^{-2}X^{-1}[Y^{-1},X^{-1}]^{N-t-1}Y^{-1}X^{-1}[X^{-1},Y^{-1}]^t
\]
after one free reduction, and this word is both reduced and cyclically reduced. To perform further reductions, we would need $u'$
to have a subword which is a cyclic  conjugate of \[[X^{-1},Y^{-1}]^tX^{-1}Y^{-1}\] or its inverse \[YX[Y^{-1},X^{-1}]^t,\] 
the first of which is not the case.

The second of these does occur as a subword of $u'$ if $N-t-1>t$. In this case we set $r_1'=r_1$ and $r_2'=r_2$ and observe that
\[YX[Y^{-1},X^{-1}]^tY=r_1'r_2'r_1'\]  is not a subword of $u'$.

{\bf Case (\ref{p:Y12}).}
In this case we have 
\[
r_1=Y^{-1},\quad r_2=[X,Y]^tX^{-1}Y^{-1}\text{ for }0\le t\le N-1.\] Then
$r_2^{-1}=YX[Y,X]^t$. Replacing $r_1r_2r_1$ by $r_2^{-1}$ gives us
\[
u'=XY[X,Y]^{N-t-1}XYX^{-1}YX[Y,X]^tX^{-1}[Y^{-1},X^{-1}]^NYX,
\]
which is both reduced and cyclically reduced. For further reductions to be possible, we would require a subword of $u'$ which is a
cyclic conjugate of \[Y^2X[Y,X]^t\] or its inverse, which  up to cyclic permutation is \[Y^{-1}[X,Y]^tX^{-1}Y^{-1}.\] Clearly this is not the case
if $t\geq 1$.

If $t=0$ then \[u'=XY[X,Y]^{N-1}XYX^{-1}YXX^{-1}[Y^{-1},X^{-1}]^NYX,\]  which reduces to
\[XY[X,Y]^{N-1}XYX^{-2}YX[Y^{-1},X^{-1}]^{N-1}YX,\] which  is both reduced and cyclically reduced. Now, $r_1r_2=Y^{-1}X^{-1}Y^{-1}$,
which does  not appear as a  subword of  $u'$ up to cyclic permutation, but  its inverse $YXY$ does. In order to perform further reductions,
we would have to either set $r_1'=Y$  and $r_1'r_2'r_1'=YXY^2$ or  $r_1'=X$  and $r_1'r_2'r_1'=XY^2X$. Neither of these words  appear,
even  after cyclic permutation.

{\bf Case (\ref{p:Y13}).}
In this case we have 
\[
r_1=Y^{-1},\quad r_2=[X,Y]^sX^{-1}Y^{-2}X^{-1}[Y^{-1},X^{-1}]^t\text{ for }0\le s,t\le N-1.\] Then
$r_2^{-1}=[X^{-1},Y^{-1}]^tXY^2X[Y,X]^s$. Replacing $r_1r_2r_1$ by $r_2^{-1}$ gives us
\[
u'=XY[X,Y]^{N-s-1}XYX^{-1}[X^{-1},Y^{-1}]^tXY^2X[Y,X]^sX^{-1}YX[Y^{-1},X^{-1}]^{N-t-1}YX,
\]
which is both reduced and cyclically reduced. To perform further  reductions, we must have a subword of $u'$ which is a cyclic
permutation of \[Y[X^{-1},Y^{-1}]^tXY^2X[Y,X]^s\] or its inverse, which is not the case.

{\bf Case (\ref{p:Y14}).}
In this case we have 
\[
r_1=Y^{-1},\quad r_2=Y^{-1}X^{-1}[Y^{-1},X^{-1}]^t\text{ for }0\le t\le N-1.\] Then
$r_2^{-1}=[X^{-1},Y^{-1}]^tXY$. Replacing $r_1r_2r_1$ by $r_2^{-1}$ gives us
\[
u'=XY[X,Y]^NX^{-1}[X^{-1},Y^{-1}]^tXYX^{-1}YX[Y^{-1},X^{-1}]^{N-t-1}YX,
\]
which is both reduced and cyclically reduced. Supposing the possibility of further reductions, we would need a subword which is a cyclic
permutation of  \[[X^{-1},Y^{-1}]^tXY^2\] or its inverse, which up to cyclic permutation is given by \[Y^{-1}X^{-1}[Y^{-1},X^{-1}]^tY^{-1}.\]
We easily see that this is not the case.
\ep

\appendix

\maketitle

\section{A computational approach, by Xinghua Gao}
As mentioned in the introduction, it was not known whether there exist points $\xi\in\fC$
such that $\rho_\xi$ has a non-free torsion-free image. In this section,
we provide and implement computational heuristics that can find such examples of
$\rho_\xi$, which turn out to have closed hyperbolic 3-manifold groups as images.

\subsection{An arithmetic formulation of the problem}
The starting point is a hyperbolic 3-manifold $T$ with a single cusp,
whose fundamental group $G$ is generated by two elements $a$ and $b$.
As an explicit example, one can take a hyperbolic $2$-bridge knot complement.
Then, consider a hyperbolic Dehn filling  $T_q$ of $T$ for some $q\in\bQ$, and a number field $\bQ(\alpha)$ such that 
\[  G=\form{a,b}=\pi_1(T_q)\le\PSL(2,\bQ(\alpha)).\]
We abuse language slightly and also write $G=\form{a,b}$ for a lift of $G$ to $\SL(2,\bQ(\alpha))$ (cf.~\cite{Kra1985,Culler1986}).

\begin{que}\label{q:computation}
Under what  conditions on $T$ and $q$ do the following conclusions hold?
\be[(i)]
\item There is a Galois automorphism $\sigma\co \alpha\mapsto\beta$ with $\beta\in\bR$, .
\item $\{\tr a^\sigma,\tr b^\sigma, \tr (ab)^\sigma\}\subset\bR\setminus[-2,2]$;
\item $\tr[a,b]^\sigma\in(-2,2)$.
\ee
\end{que}

Computationally, Question~\ref{q:computation} suggests that we enumerate such possible $T$ and $q$,
and verify all of the above three conditions. 
As we have the hyperbolic structure of $T_q$ and the algebraic number $\alpha\in\bC$,
the verification step should be computationally straightforward up to the computation of Galois conjugates. 
Once we have such an example of $T$ and $q$, then the resulting point
\[
\xi:=(\tr a^\sigma,\tr b^\sigma,\tr(ab)^\sigma)
\]
is a point in the character variety of $\fC_\theta$ for $\theta$ satisfying
\[-2\cos(\theta/2) = \tr[a,b]^\sigma.\]
Moreover, the image of the monodromy $\rho_\xi$ is isomorphic to $\pi_1(T_q)$ as we desire.

\begin{rem}
Note that a parabolic generator of $\pi_1(T)$ will remain parabolic after Galois conjugation.
\end{rem}

\subsection{An algorithm for producing explicit examples}
Examples of $(T,q)$ with  the desired properties can be produced using SnapPy~\cite{SnapPy,sage}.
In order to construct a representation with the image in $\PSL(2,\mathbb{R})$, we will need the following fact.

\begin{lem}\cite[Corollary 3.2.5]{MR2003book}
\label{l:psl2r}
If $\Gamma$ is a nonelementary subgroup of $\SL(2,\mathbb{C})$ such that
$\mathbb{Q}(\text{tr } \Gamma)$ is a subset of $\mathbb{R}$, then $\Gamma$ is conjugate to a subgroup of $\SL(2,\mathbb{R})$.
\end{lem}

Let $T_r$ be a closed hyperbolic $3$-manifold obtained by applying Dehn filling
to one-cusped hyperbolic $3$-manifold $T$ along a curve of slope $r$. 
We fix a triangulation of $T$ by ideal tetrahedra, which also gives an ideal triangulation of  $T_r$.
Then the edge gluing equations of $T$ together with the Dehn filling equation determine the hyperbolic
structure of $T_r$. The solution to the this system of equations is a set of complex numbers,
which parameterize the shape of tetrahedra in $T_r$. These parameters generate the
\emph{tetrahedral field} $K$ of $T_r$. We want to find a real embedding $\sigma$ of the number
field $K$ so that all the tetrahedra of $T$ have real shapes after applying $\sigma$. 
The associated holonomy representation then gives a faithful representation
\[\rho_{\mathbb{R}}\colon \pi_1(T_r)\to\PSL(2,\mathbb{C})\] with all matrices in the image having real trace.
Therefore by Lemma~\ref{l:psl2r},  the group $\rho_{\mathbb{R}}(\pi_1(T_r))$ is conjugate into $\PSL(2,\mathbb{R})$.
Since a conjugacy leaves the trace unchanged, we obtain a desired hyperbolic $3$-manifold $T_r$ giving an affirmative answer to
Question~\ref{q:computation}.

 The process of finding a suitable $T$ and corresponding tetrahedral field can be formulated with the following steps:

\begin{itemize}
\item[Step 1)] Let $T$ be a hyperbolic knot complement with the default triangulation in SnapPy and apply
 Dehn filling of slope $r$. Compute the tetrahedral field $K$ of $T_r$ using the SnapPy manifold class
 \texttt{tetrahedra\_field\_gens()}. We can then use the SageMath number field class
 \texttt{find\_field()} to find the defining polynomial of $K$. \\  
\item[Step 2)] Find a real embedding $\sigma$ of $K$, if there exits one, using the SageMath number field class \texttt{real\_embeddings()}.\\
\item[Step 3)] Apply the real embedding $\sigma$ and set up the new triangulation with
real shape parameters, using the SnapPy manifold class \texttt{set\_tetrahedra\_shapes()}.\\
\item[Step 4)] Computes the associated holonomy representation $\rho_{\mathbb{R}}$
and the image of $a$, $b$, $ab$ and $[a,b]$ under $\rho_{\mathbb{R}}$, using the SnapPy fundamental group classe \texttt{SL2C()}.
Finally, compute the resulting traces.\\ 
\end{itemize}

\subsection{Explicit examples}
Here we produce several examples giving affirmative answer to Question~\ref{q:computation} .
The traces are truncated to four places after the decimal point.

\begin{table}[H]
\begin{tabular}{@{}llllll@{}}
\hline
\multicolumn{1}{r}{}
 & $\text{tr} \rho_{\mathbb{R}}(a)$ & $\text{tr} \rho_{\mathbb{R}}(b)$ & $\text{tr} \rho_{\mathbb{R}}(ab)$  & $\text{tr} \rho_{\mathbb{R}}([a,b])$ \\
\hline
$7_6(0)$        & 2.4509  & 2.0881  & 2.4509  &1.8307\\
$8_{13}(0)$    & 2.1258  & 2.7610  & 2.4523  & 1.7623  \\
$9_{12}(0)$    & 2.0382  & -2.4497 & -2.4497 & 1.9249  \\
$9_{15}(0)$    & -2.2535 & 2.1399  & -2.2535 & 1.8686 \\
$10_{10}(0)$  & -3.7588 & -3.0575 & 9.0343  &  -0.7349 \\
\hline
\end{tabular}
\end{table}

In this table, $a$ and $b$ are the generators of the corresponding fundamental group, with the default triangulation in SnapPy.
Presentations of the fundamental group of $7_6(0)$ and $8_{13}(0)$ are included below, as well as the matrix representatives of
$a$ and $b$. Presentations of the fundamental groups of the other three manifolds are unwieldy due to their size, so we have omitted them.
The interested readers may use the manifold class  \texttt{fundamental\_group()}  to verify our claims.

For the first manifold,
\[\begin{split}
\pi_1(7_6(0))=\left< a,b \right. | &abABBAbABabbaBabbaBAbABBAbaBabbaB,  \\
  & \left. aaaaabABBAbABabbaBabAbaBabbaBAbABBAb \right> ,
\end{split}\]

\begin{displaymath}
\rho_{\mathbb{R}}(a)=
\begin{pmatrix}
0.5171 & 0 \\
-0.3455 & 1.9338 
\end{pmatrix}, \quad
\rho_{\mathbb{R}}(b)=
\begin{pmatrix}
1.0881 & 0.1319 \\
0.6682 & 1 
\end{pmatrix}.
\end{displaymath}
Here, upper case and  lower case versions of a letter are inverses of each other.
Again, we truncate matrix entries to four places.
We note however that they are all algebraic numbers in a real embedding of the tetrahedra field $K$ of $7_6(0)$,
with the defining polynomial
\[\begin{split} p(x)=x^{12} - 2x^{11} - 2x^{10} + 14x^9 - 25x^8 + 32x^7 \\- 35x^6 + 38x^5 - 38x^4 + 30x^3 - 17x^2 + 6x - 1.\end{split}\] 

For the second manifold,
\[\begin{split}
\pi_1(8_{13}&(0))=\\
\left< a,b \right. | &aaabABBAAABabbaaabABBAbaaabbaBAAABBAbaaabbaBAAABabb,  \\
  & \left. aaabABBAAABabbaaabABBAbABBAbaaabbaBAAABBAbaaabABAb \right> ,
\end{split}\]

\begin{displaymath}
\rho_{\mathbb{R}}(a)=
\begin{pmatrix}
1.1258 & 0.3547i \\
-0.3547i & 1 
\end{pmatrix}, \quad
\rho_{\mathbb{R}}(b)=
\begin{pmatrix}
2.8986 & -2.4657i \\
-0.5673i & -0.1376
\end{pmatrix}.
\end{displaymath}
Note that these two matrices are not in $\PSL(2, \mathbb{R})$, but they are simultaneously conjugate into
$\PSL(2, \mathbb{R})$, according to Lemma \ref{l:psl2r}. 
The defining polynomial of the tetrahedra field $K$ of $8_{13}(0)$ is
\[\begin{split}
p(x)=x^{14} - x^{13} - 5x^{12} - 4x^{11} + 10x^{10} + 14x^9 - 10x^8\\ - 29x^7 - 5x^6 + 29x^5 + 19x^4 - 11x^3 - 17x^2 - 7x - 1.\end{split}\]

For now, we are only able to produce closed hyperbolic $3$-manifolds via $0$-Dehn filling to two-bridge knot complements as examples.

\section*{Acknowledgements}
The authors thank A.~Reid and R.~Schwartz for valuable discussions. 
X.~Gao would like to thank K.~Ohshika for suggesting to her the question in the appendix and thank N.~Dunfield for helping her with the  code. 
S.~Kim and X.~Gao are supported by Mid-Career Researcher Program (2018R1A2B6004003) through the
National Research Foundation funded by the government of Korea.  
S.~Kim is supported by Samsung Science and Technology Foundation under Project Number SSTF-BA1301-51, and by KIAS Individual Grant (MG073601) at Korea Institute for Advanced Study.
T.~Koberda is partially supported  by an Alfred P. Sloan Foundation Research Fellowship and by NSF Grant DMS-1711488.
J.~Lee is supported by the grant NRF-2019R1F1A1047703.
K.~Ohshika is partially supported by the JSPS Grant-in-aid for Scientific Research No 17H02843.
Tan is partially supported by the National University
of Singapore academic research grant R-146-000-289-114.
The authors thank the anonymous referee for helpful suggestions.

\bibliographystyle{amsplain}
\bibliography{ref}\end{document}